\input ./style/arxiv-general.cfg
\documentclass[aap,MSNbibl,dvips]{arximspdf}
\makeatletter
   \@ifpackageloaded{graphicx}{}{\usepackage{graphicx}}
\makeatother
\doi{10.1214/14-AAP1093}% Updated by VTEXPTS2LaTeX.exe, 06.02.2015
\volume{26}
\issue{1}
\pubyear{2016}
\firstpage{328}
\lastpage{359}
\docsubty{FLA}

\makeatletter
\newcommand{\tr}{\mathrm{tr}}
\newcommand{\lleft}{\left}
\newcommand{\rrvert}{\vert}
\newcommand{\rright}{\right}
\newcommand{\llvert}{\vert}
\def\ind{{\mathbf1}}
\def\re{\operatorname{Re}}
\def\im{\operatorname{Im}}
\def\p{{\mathbb P}}
\def\e{{\mathbb E}}
\def\r{{\mathbb R}}
\def\c{{\mathbb C}}
\def\d{\mathrm{d}}
\def\i{\mathrm{i}}
\newtheorem{theorem}{Theorem}
\newtheorem{lemma}{Lemma}
\newtheorem{proposition}{Proposition}
\newproclaim{definition}{Definition}
\newproclaim{ass}{Assumption}
\newproclaim{rem}{Remark}
\makeatother

\begin{document}
\begin{frontmatter}

\title{Approximating L\'{e}vy processes with completely monotone jumps}
\runtitle{Approximating L\'{e}vy processes}

\begin{aug}
% Corresponding author: Alexey Kuznetsov - kuznetsov@mathstat.yorku.ca% Updated by VTEXPTS2LaTeX.exe, 12.02.2015 17:23
%Updated by VTEXPTS2LaTeX.exe, 06.02.2015 09:09
\author[A]{\fnms{Daniel}~\snm{Hackmann}\ead[label=e1]{dhackman@mathstat.yorku.ca}\ead[label=u1,url]{www.danhackmann.com}\thanksref{T1}}
\and
\author[A]{\fnms{Alexey}~\snm{Kuznetsov}\corref{}\ead[label=e2]{kuznetsov@mathstat.yorku.ca}\ead[label=u2,url]{www.math.yorku.ca/\texttildelow akuznets}\thanksref{T2}}
\runauthor{D. Hackmann and A. Kuznetsov}
\affiliation{York University}
%\dedicated{}
\address[A]{Department of Mathematics and Statistics\\
York University\\
4700 Keele Street\\
Toronto, ON, M3J 1P3\\
Canada\\
\printead{e1}\\
\phantom{E-mail: }\printead*{e2}\\
\printead{u1}\\
\phantom{URL: }\printead*{u2}}
\end{aug}
\thankstext{T1}{Supported by the Ontario Graduate Scholarship Program.}
\thankstext{T2}{Supported by the Natural Sciences and Engineering Research Council of Canada.}

% HISTORY:
%
\received{\smonth{4} \syear{2014}}% Updated by VTEXPTS2LaTeX.exe,
%06.02.2015 09:09
%
\revised{\smonth{12} \syear{2014}}% Updated by VTEXPTS2LaTeX.exe,
%06.02.2015 09:09

% ABSTRACT
%
\begin{abstract}
L\'evy processes with completely monotone jumps appear frequently in various
applications of probability. For example, all popular stock price
models based on L\'evy processes
(such as the Variance Gamma, CGMY/KoBoL and Normal Inverse Gaussian)
belong to this class.
In this paper we continue the work started in
[\textit{Int. J. Theor. Appl. Finance} \textbf{13} (2010) 63--91,
\textit{Quant. Finance} \textbf{10} (2010) 629--644]
and develop a simple
yet very efficient method for approximating processes with completely
monotone jumps by
processes with hyperexponential jumps, the latter being the most
convenient class for
performing numerical computations. Our approach is based on connecting
L\'evy processes
with completely monotone jumps with several areas of classical
analysis, including
Pad\'e approximations, Gaussian quadrature and orthogonal polynomials.
\end{abstract}

% KEYWORDS
\begin{keyword}[class=AMS]
\kwd[Primary ]{60G51}
%\kwd{}
\kwd[; secondary ]{26C15}
\end{keyword}
\begin{keyword}
\kwd{L\'evy processes}
\kwd{complete monotonicity}
\kwd{hyperexponential processes}
\kwd{Pad\'e approximation}
\kwd{rational interpolation}
\kwd{Gaussian quadrature}
\kwd{Stieltjes functions}
\kwd{Jacobi polynomials}
\end{keyword}
\end{frontmatter}

\setcounter{footnote}{2}
%s1 #&#
\section{Introduction}\label{secintroduction}\label{sec1}

Most researchers working in Applied Mathematics are familiar with the
problem of choosing the right mathematical objects for their modeling
purposes: one needs to strike a balance between the simplicity of the
model, its analytical and numerical tractability and its ability to
provide a realistic \mbox{description} of the phenomenon. For example, when
modeling stock prices in mathematical finance we are faced with the
following dilemma: do we choose a process which fits the empirically
observed behavior of stock prices (such as having jumps of infinite
activity \cite{carr2002fine}), or do we settle for a simpler model
which provides for explicit formulas and efficient numerical algorithms?
The first choice would lead to the most popular families of L\'evy
processes, such as the Variance Gamma (VG), CGMY/KoBoL, Meixner and
Normal Inverse Gaussian (NIG) families. These processes, which belong
to a wider class of processes with completely monotone jumps, provide a
good fit for market data, and they are flexible enough to accommodate
for such desirable features as jumps of infinite activity and finite or
infinite variation. They also enjoy a certain degree of analytical
tractability (e.g., VG and NIG processes have explicit transition
probability densities), and European option prices and Greeks can be
computed quite easily. However, the computation of more exotic option
prices (such as barrier, lookback and Asian options) is a much more
challenging task. On the other hand, hyperexponential processes (also
known as ``hyperexponential jump-diffusion processes,'' see \cite
{CaiKou}), and more general processes with jumps of rational transform
(see \cite{Fourati,Kuznetsov,Mordecki}) form the most convenient class
for performing numerical calculations. This
is due to the fact that these processes have an explicit Wiener--Hopf
factorization, which leads to simple and efficient numerical algorithms
for pricing barrier and lookback options \mbox{\cite{CaiKou,Jeannin}} and
Asian options \cite{CaiKouAsianoptions}. One might think that
hyperexponential processes are perfect candidates for modeling stock
prices, yet they have a major flaw in that their jumps are necessarily
of finite activity, which seems to be incompatible with empirical
results \cite{carr2002fine}.

A natural way to reconcile these two competing objectives is to
approximate processes with completely monotone jumps with the
hyperexponential processes. Two approximations of this sort were
developed recently: Jeannin and Pistorius \cite{Jeannin}
use the least squares optimization in order to find the approximating
hyperexponential process, while
Crosby, Le Saux and Mijatovi\'c \cite{Crosby} use a more direct
approach based on the Gaussian quadrature. Our goal in this paper is to
present a new method for approximating L\'evy processes with completely
monotone jumps, and to demonstrate that this method is natural, simple
and very efficient.

Let us present the main ideas behind our approach. Approximating a L\'
evy process $X$ is equivalent to finding an approximation to
its \textit{Laplace exponent}, defined as $\psi(z):=\ln\e[ \exp
(zX_1)]$. The Laplace exponent of a hyperexponential processes is a
rational function, therefore, our problem reduces to two steps: (i)
finding a good rational approximation $\tilde\psi(z) \approx\psi
(z)$, and (ii) ensuring that the rational function $\tilde\psi(z)$ is
itself a Laplace exponent of some L\'{e}vy process $\tilde X$. For the
first step, we rely on the extensive literature on rational
approximations and interpolations. One of the simplest and the most
natural methods of rational approximation is the \textit{Pad\'e
approximation}; see the classical book by Baker \cite{Baker} for an
excellent account of this theory. The Pad\'e approximation
$f^{[m/n]}(x)$, of a function $f(x)=\sum_{k\ge0} c_n x^k$, is defined
as a rational function $P_m(x)/Q_n(x)$ [where $P$ and $Q$ are
polynomials satisfying $\deg(P)\le n$ and $\deg(Q) \le m$] which
matches the first $n+m+1$ Taylor coefficients of $f(x)$. Pad\'e
approximations are easy to compute, and there exists a well developed
theory related to their various properties (convergence, error
estimates, etc.). Thus, the first step of our program is rather simple,
but the second step is much more challenging: we need to ensure that
the approximating rational function $\psi^{[m/n]}(z)$ is itself the
Laplace exponent of some
L\'evy process $Y$. First, we can considerably reduce the number of
possible cases that we need to study. It is known
(see Proposition 2 on page 16 in \cite{Bertoin})
that the Laplace exponent of a L\'evy process satisfies $\psi(\i z) =
O(z^2)$ as $z\to\infty$, therefore, the functions
$\psi^{[m/n]}(z)$ cannot be Laplace exponents if $m>n+2$. If $m<n$
then necessarily $\psi^{[m/n]}(\i z) \to0$ as $z\to\infty$, and
one can
prove\footnote{Assume that a L\'evy process $Y$ has a rational
function $f(z)=P(z)/Q(z)$ as its Laplace exponent. Do the partial
fraction decomposition of $f(z)$ and identify the L\'evy measure of $Y$
via the L\'evy--Khintchine formula and the inverse Laplace transform.
Show that if $\lim_{z \to\infty} f(z)=\lambda<\infty$, then $Y$
must be a compound Poisson process with jump intensity $\lambda$. In
particular, if $\lim_{z \to\infty} f(z)=0$ then $Y= 0$ almost surely.}
that a rational function with this property cannot be the Laplace
exponent of a L\'evy process $Y$ (unless $Y = 0$ almost surely).
This shows that in the full table of Pad\'e approximations $\psi
^{[m/n]}(z)$ only the following functions:
%
%e1 #&#
%
\begin{equation}
\label{threePadeapproximations} \psi^{[n/n]}(z),\qquad \psi^{[n+1/n]}(z)\quad\mbox{and}\quad
\psi^{[n+2/n]}(z)
\end{equation}
can possibly be Laplace exponents of a L\'evy process.

Checking whether a given function is the Laplace exponent of a L\'evy
process is a very difficult task: one would need to show that the
function can be represented via the L\'evy--Khintchine formula [see
formula (\ref{LevyKhinchine}) below]. Since it is impossible to
verify this property numerically, one would require some additional
qualitative information about the function. In our case, this
additional information comes from the fact that $\psi(z)$ is the
Laplace exponent of a process with \textit{completely monotone jumps}.
Using this key fact
and utilizing connections with several branches of classical analysis
(such as the theory of Pad\'e approximations, orthogonal polynomials,
Stieltjes functions and Gaussian quadrature), we are able to completely
characterize all cases when the functions in
(\ref{threePadeapproximations}) are Laplace exponents of L\'evy
processes. Our main result states that if the original L\'evy process
has completely monotone jumps, the function $\psi^{[n+1/n]}(z)$ is a
Laplace exponent of a hyperexponential L\'evy process $X^{(n)}$, which
converges to $X$ in distribution as $n\to+\infty$. Moreover, if the
process $X$ has only positive (or only negative) jumps, the same
results holds true for $\psi^{[n+2/n]}(z)$ [and for $\psi^{[n/n]}(z)$
under the additional assumption that the process has jumps of finite variation].

The paper is organized as follows. Section~\ref{sectionmainresults}
contains our main results on approximating L\'evy processes with
completely monotone jumps (treating the two-sided and one-sided cases
separately). Section~\ref{sectionexplicitexamples} discusses the
important special cases of the Gamma subordinator and of the one-sided
tempered stable processes; in both cases the Pad\'e approximation is
given explicitly. In this section we also discuss how to use these
results to construct explicit approximations to VG, CGMY and NIG
processes, and present some extensions of our approximation scheme,
including (i) the use of Pad\'e approximation centered at an arbitrary
point and (ii) a more general multi-point rational interpolation
technique. In Section~\ref{sectionnumerics}, we present the results
of several numerical experiments which demonstrate the efficiency of
our approximation method. We compute the L\'evy density, the CDF and
the prices of various options for the approximating processes and
investigate their convergence. In Section~\ref{sectionremarks}, we
compare our approach with the methods developed in \cite{Jeannin} and
\cite{Crosby} and we discuss connections with meromorphic processes.
For the reader's convenience, in the \hyperref[appendix1]{Appendix} we collect
some results from the theory of Pad\'e approximations, Stieltjes
functions, Gaussian quadrature and orthogonal polynomials, which are
used elsewhere in this paper.

%s2 #&#
\section{Main results}\label{sectionmainresults}

We begin by introducing a number of key definitions and notations. Let
$X$ be a L\'evy process, and let $\psi(z):= \ln\e[e^{zX_1}]$ denote
its \textit{Laplace exponent}, which is initially defined on the vertical
line $z\in\c$, $\re(z)=0$. The L\'evy--Khintchine formula states that
%
%e2 #&#
%
\begin{equation}
\label{LevyKhinchine} \psi(z)=\sigma^2z^2/2+ a z + \int
_{\r} \bigl(e^{zx}-1-zh(x) \bigr) \Pi(\d x),
\end{equation}
where $\sigma\ge0$, $a\in\r$, the L\'evy measure $\Pi(\d x)$
satisfies $\int_{\r} (1 \wedge x^2) \Pi(\d x)<\infty$, and $h(x)$
is the \textit{cutoff function}, which is required to ensure the
convergence of the integral. Everywhere in this paper we will work
under the following assumption.

\begin{ass}\label{ass1}
The L\'evy measure $\Pi(\d x)$ is absolutely continuous, and its
density $\pi(x)$ decreases exponentially fast as $x \to\pm\infty$.
\end{ass}

If the cutoff function is fixed [the classical choice is $h(x)\equiv x
\ind_{\{\llvert  x\rrvert  <1\}}$], then the process $X$ is completely characterized
by the triple $(a,\sigma^2, \pi)$, which determines the Laplace
exponent in (\ref{LevyKhinchine}). It is often convenient, however,
to use different cutoff functions depending on the situation.
Everywhere in this paper we will follow the convention that if the
process $X$ has jumps of finite variation, we will take $h(x)\equiv0$,
otherwise we will set $h(x) \equiv x$ (which is a legitimate choice due
to Assumption~\ref{ass1}). To distinguish between these two cases we will write the
characteristic triple as $(a,\sigma^2,\pi)_{h\equiv0}$ in the former
case and $(a,\sigma^2,\pi)_{h \equiv x}$ in the latter case.

We recall that a function $f\dvtx  (0,\infty) \mapsto\r$ is called \textit{completely monotone} if $(-1)^k f^{(k)}(x)\ge0$ for all
$k=0,1,2,\dots$ and $x>0$.

%de1 #&#
%
\begin{definition}\label{defcompletelymonotone}
We say that the process $X$ has \emph{completely monotone jumps} if
the functions $\pi(x)$ and $\pi(-x)$ are completely monotone
for $x \in(0,\infty)$.
\end{definition}

Using Bernstein's theorem (see \cite{SSV2012}, page 3), we can express
the above condition in an equivalent form:
$X$ has completely monotone jumps if and only if there exists a
positive Radon measure $\mu$, with support in $\r\setminus\{0\}$,
such that for all $x\in\r$
%
%e3 #&#
%
\begin{equation}
\label{defpi} \pi(x)=\ind_{\{x>0\}} \int_{(0,\infty)}
e^{-u x}\mu(\d u) + \ind _{\{x<0\}} \int_{(-\infty,0)}
e^{-u x} \mu(\d u).
\end{equation}
For our further results, we will need the following two facts [which
follow easily from (\ref{defpi}) by Fubini's theorem]:
%
%e4 #&#
%e5 #&#
%
\begin{eqnarray}
\label{generalcondition}
\int_{\r} x^2 \pi(x)\, \d x &<& \infty \quad\mbox{if and only if}\quad \int_{\r} \llvert u\rrvert ^{-3} \mu(\d u)<\infty,
\\
\label{finvariationcondition}
\int _{\r} \llvert x\rrvert \pi(x) \,\d x&<& \infty\quad\mbox{if and only if}\quad \int _{\r} u^{-2} \mu(\d u)<\infty.
\end{eqnarray}
Condition (\ref{generalcondition}) is required to ensure that the
function $\pi(x)$ can be considered as a L\'evy density, while the
stronger condition (\ref{finvariationcondition}) ensures that the
resulting L\'evy process $X$ has jumps of finite variation.

Assuming that the L\'evy density $\pi(x)$ is given by (\ref{defpi}),
we denote
\begin{eqnarray*}
\rho &:=& \sup \biggl\{c \ge0\dvtx  \int_{\r^+} e^{cx}
\pi(x) \,\d x < \infty \biggr\}=\sup \bigl\{u\ge0\dvtx  \mu\bigl((0,u)\bigr)=0 \bigr\},
\\
\hat\rho &:=& \sup \biggl\{c \ge0\dvtx  \int_{\r^-} e^{-cx}
\pi(x) \,\d x < \infty \biggr\}=\sup \bigl\{u\ge0\dvtx  \mu\bigl((-u,0)\bigr)=0 \bigr\}.
\end{eqnarray*}
Assumption \ref{ass1} implies that $\rho>0$ and $\hat\rho>0$. We will denote
by ${\mathcal{CM}}(\hat\rho,\rho)$ the class of L\'evy processes
with completely monotone jumps and parameters $\rho$ and $\hat\rho$
defined as above.

Now we consider an important subclass of ${\mathcal{CM}}(\hat\rho,\rho)$.
%de2 #&#

\begin{definition}\label{defhyperexponential}
We say that the process $X$ has \emph{hyperexponential jumps} if the
support of the
measure $\mu(\d x)$ in (\ref{defpi}) consists of finitely many points.
\end{definition}

Let us consider a hyperexponential process $X$. According to Definition
\ref{defhyperexponential}, the measure $\mu$
has finite support,\vspace*{1pt} which we will denote $\operatorname{supp}(\mu)=\{
\hat\beta_i\}_{1\le i \le\hat N} \cup
\{\beta_i\}_{1\le i \le N}$, where $\hat N\ge0$ and $N \ge0$, and
$\hat\beta_i<0$ and $\beta_i>0$. We denote $\mu(\{\hat\beta_i\}
)=\hat\alpha_i$ and $\mu(\{\beta_i\})=\alpha_i$. Then
the L\'evy density of $X$ can be represented in the form
%
%e6 #&#
%
\begin{equation}
\label{defhyperexponentialpi} \pi(x)={\mathbf1}_{\{x>0\}} \sum
_{i=1}^N
\alpha_i e^{-\beta_i
x}+{\mathbf1}_{\{x<0\}} \sum
_{i=1}^{\hat N} \hat\alpha_i
e^{-\hat\beta_i x},
\end{equation}
where one of the sums can be empty (if $\hat N=0$ or $N=0$). Formula
(\ref{defhyperexponentialpi}) provides another equivalent definition
of a hyperexponential process, as having positive/negative jumps equal
in law to a finite mixture of exponential distributions.

%de3 #&#
%
\begin{definition}\label{defpade}
Let $f$ be a function with a power series representation $f(z) = \sum_{i=0}^{\infty}c_i(z-a)^i$.
If there exist polynomials $P_m(z)$ and $Q_n(z)$ satisfying $\deg
(P)\le m$, $\deg(Q)\le n$, $Q_n(a) \ne0$
and
\begin{eqnarray*}
&&\frac{P_m(z)}{Q_n(z)} = f(z) + O\bigl((z-a)^{m+n+1}\bigr),\qquad z\to a,
\end{eqnarray*}
then we say that $f^{[m/n]}(z):= P_m(z)/Q_n(z)$ is the \emph{$[m/n]$
Pad\'{e} approximant of $f$ at point $a$}.
\end{definition}
Everywhere in this paper we will consider the case when the power
series representation for $f(z)$ is convergent in some neighborhood of
$a$ (more\vspace*{1pt} generally, it can also be considered as a formal power
series). When $a=0$, we will call $f^{[m/n]}(z)$ simply the
$[m/n]$ Pad\'e approximation of $f(z)$, without mentioning the
reference point.

%s2.1 #&#
\subsection{Approximating L\'evy processes with two-sided jumps}

For a L\'evy process $X \in{\mathcal{CM}}(\hat\rho,\rho)$, we define
%
%e7 #&#
%
\begin{equation}
\label{defmustar} \mu^*(A)=\mu\bigl( \bigl\{ v \in\r\dvtx  v^{-1} \in A\bigr\}
\bigr),
\end{equation}
for all Borel sets $A\subset\r$, where the measure $\mu(\d v)$
appears in (\ref{defpi}).
Note that $\operatorname{supp}(\mu^*)\subseteq[-1/\hat\rho, 1/\rho]$,
and if the measure $\mu(\d v)$ is absolutely continuous with a density $m(v)$,
then $\mu^*(\d v)$ also has a density, given by
$m^*(v)= m(1/v)/v^2$.
The measure $\mu^*(\d v)$ will play a very important role in this paper.
%le1 #&#

\begin{lemma}\label{lemmamustar}
Assume that $X \in{\mathcal{CM}}(\hat\rho, \rho)$. Then
\begin{eqnarray*}
&&\int
_{[-1/\hat\rho,1/\rho]} \llvert v\rrvert ^3 \mu^*(\d v)<\infty,
\end{eqnarray*}
and
\[
\int
_{[-1/\hat\rho,1/\rho]} v^2 \mu^*(\d v)<\infty \quad\mbox{if and only if}\quad  X\mbox{ has jumps of finite variation}.
\]
\end{lemma}

\begin{pf}
The result follows from (\ref{generalcondition}) and (\ref
{finvariationcondition}) by change of variables $u=1/v$.
\end{pf}

Now we are ready to introduce our first approximation. We start with a
L\'evy process
$X \in{\mathcal{CM}}(\hat\rho, \rho)$ defined by the
characteristic triple
$(a,0,\pi)_{h\equiv x}$. Note that the process $X$ has zero Gaussian
component. However, there is no lack of generality in assuming this: if
we know how to approximate a L\'evy processes with zero Gaussian
component, we know how to approximate a general L\'evy process, as we
can always add a scaled Brownian motion to our hyperexponential approximation.

According to Lemma \ref{lemmamustar}, $\llvert  v\rrvert  ^3\mu^*(\d v)$ is a
finite measure on the interval $[-1/\hat\rho,  1/\rho]$. Let $\{x_i\}
_{1\le i \le n}$ and $\{w_i\}_{1\le i \le n}$ be the nodes and the
weights of the Gaussian quadrature of order $n$ with respect to this
measure (we have included the definition and several key properties of
the Gaussian quadrature in the \hyperref[appendix1]{Appendix}.
We define
%
%e8 #&#
%
\begin{equation}
\label{defpsin} \psi_{n}(z):=az + z^2 \sum_{i=1}^n \frac{w_i}{1-zx_i}.
\end{equation}

%th1 #&#
%
\begin{theorem}\label{thmmain}
\textup{(i)} The function $\psi_{n}(z)$ is the $[n+1/n]$ Pad\'e
approximant of~$\psi(z)$.

\textup{(ii)} The function $\psi_n(z)$ is the Laplace exponent of a
hyperexponential process $X^{(n)}$ with the characteristic triple
$(a,\sigma_n^2,\pi_n)_{h\equiv x}$, where
%
%e9 #&#
%
\begin{eqnarray}
\label{defsigman} \sigma_n^2:= \cases{ 0, &\quad if
$x_i \ne0$ for all $1\le i \le n$,
\cr
2w_j, &\quad if
$x_j = 0$ for some $1\le j \le n$,}
\end{eqnarray}
and
%
%e10 #&#
%
\begin{eqnarray}
\label{defpin} \pi_n(x):= \cases{ \displaystyle\sum
_{1\le i \le n\dvtx  x_i<0} w_i \llvert x_i\rrvert
^{-3} e^{-x/x_i}, &\quad if $ x<0$,
\vspace*{6pt}\cr
\displaystyle\sum
_{1\le i \le n\dvtx  x_i>0} w_i x_i^{-3}
e^{-x/x_i}, &\quad if $x>0$.}
\end{eqnarray}
If one of the sums in (\ref{defpin}) is empty, it should be
interpreted as zero.\vspace*{1pt}

\textup{(iii)} The random variables $X^{(n)}_1$ and $X_1$ satisfy $\e
[(X^{(n)}_1)^j]=\e[(X_1)^j]$ for $1\le j \le2n+1$.
\end{theorem}

\begin{pf}
Our first goal is to establish an integral representation of $\psi(z)$
in terms of the measure $\mu^*(\d v)$.
Assume that $z\in\c$ with $-\hat\rho<\re(z)<\rho$. We substitute
(\ref{defpi}) into (\ref{LevyKhinchine}), use Fubini's theorem to
interchange the order of integration and obtain
%
%e11 #&#
%
\begin{equation}
\label{formulapsiStieltjesfunction0} \psi(z)= a z + z^2 \int_{\r}
\frac{{\operatorname{sign}}(u)}{u-z} \frac{\mu(\d u)}{u^2}.
\end{equation}
Changing the variable $v=u^{-1}$ in the above integral and using the
fact that $\mu((-\hat\rho, \rho))=0$, we obtain
%
%e12 #&#
%
\begin{equation}
\label{formulapsiStieltjesfunction} \psi(z)= a z + z^2 \int
_{[-1/\hat\rho,1/\rho]} \frac
{\llvert  v\rrvert  ^3 \mu^*(\d v)}{1-vz},\qquad  -\hat\rho<\re(z)<\rho.
\end{equation}
By analytic continuation, we can see that the above formula is valid in
a larger region $\c\setminus\{(-\infty,-\hat\rho] \cup[\rho,
\infty) \}$.

Let us prove (i). By definition, the Gaussian quadrature of order $n$
is exact for polynomials of degree not greater than $2n-1$,
therefore,
\[
\int
_{[-1/\hat\rho,1/\rho]} v^k \llvert v\rrvert ^3 \mu^*(\d
v)= \sum
_{i=1}^n x_i^k
w_i,\qquad k=0,1,2,\dots, 2n-1.
\]
The above identity is equivalent to
%
%e13 #&#
%
\begin{equation}
\label{equalitymoments} \Biggl(\frac{\d^k}{\d z^k} \int
_{[-1/\hat\rho,1/\rho]} \frac{\llvert  v\rrvert  ^3 \mu^*(\d v)}{1-vz} \Biggr)
\Biggl\llvert _{z=0}= \Biggl(\frac{\d^k}{\d z^k} \sum
 _{i=1}^n
\frac{w_i}{1-zx_i} \Biggr) \Biggr\rrvert _{z=0},
\end{equation}
for $k=0,1,2,\dots,2n-1$. Formulas (\ref{defpsin}), (\ref
{formulapsiStieltjesfunction}) and (\ref{equalitymoments}) imply that
%
%e14 #&#
%
\begin{equation}
\label{psipsin} \psi^{(k)}(0)=\psi_n^{(k)}(0),\qquad
k=0,1,2,\dots, 2n+1.
\end{equation}
By definition (\ref{defpsin}), $\psi_n(z)$ is
a rational function, which can be written in the form $P(z)/Q(z)$ with
$\deg(P)\le n+1$ and $\deg(Q)=n$. Using this fact and
formula~(\ref{psipsin}), we see that $\psi_n(z)\equiv\psi
^{[n+1/n]}(z)$, which proves (i).

From (\ref{LevyKhinchine}) and (\ref{defhyperexponentialpi}), we
see that the Laplace exponent of a
hyperexponential process $Y$ having triple $(a,0,\pi)_{h\equiv x}$ is
given by
\begin{eqnarray*}
&&\psi_Y(z)=az + z^2 \sum
_{i=1}^{\hat N}
\frac{\hat\alpha
_i}{\llvert  \hat\beta_i\rrvert  ^3(1-z/\hat\beta_i)} + z^2 \sum
_{i=1}^{N}
\frac{\alpha_i}{\beta_i^3(1-z/\beta_i)}.
\end{eqnarray*}
The result of item (ii) follows at once by comparing the above
expression with (\ref{defpsin}).

Now that we have established that $\psi_n(z)$ is the Laplace exponent
of a hyperexponential process $X^{(n)}$,
formula (\ref{psipsin}) shows that the first $2n+1$ cumulants of
$X^{(n)}$ are equal to the corresponding cumulants of $X_1$, which is
equivalent to the equality of corresponding moments and proves item (iii).
\end{pf}

The next important question that we need to address is how fast the
approximations $\psi_n(z)$ converge to $\psi(z)$.
As we have seen in the proof of Theorem \ref{thmmain} (see also \cite
{Kwasnicki,Rogers}), the Laplace exponent $\psi(z)$ of a process $X
\in{\mathcal{CM}}(\hat\rho, \rho)$ is analytic in the cut complex
plane $\c\setminus\{(-\infty,-\hat\rho] \cup[\rho, \infty) \}$.
As we will establish in the next theorem, $\psi_n(z)$ converge to
$\psi(z)$ \textit{everywhere} in this region, and the convergence is
exponentially fast on compact subsets of $\c\setminus\{(-\infty,-\hat\rho] \cup[\rho, \infty) \}$. This behavior
should be compared with Taylor approximations, which can converge only
in a circle of finite radius [lying entirely in the region of
analyticity of $\psi(z)$]. This demonstrates that Pad\'e
approximations are very well suited to approximate Laplace exponents of
processes in ${\mathcal{CM}}(\hat\rho, \rho)$.

%th2 #&#
%
\begin{theorem}\label{theoremconvergence}
For any compact set $A \subset\c\setminus\{ (-\infty, -\hat\rho]
\cup[\rho, \infty)\}$, there exist $c_1=c_1(A)>0$ and $c_2=c_2(A)>0$
such that for all $z \in A$ and all $n\ge1$
\begin{eqnarray*}
&&\bigl\llvert \psi_n(z)-\psi(z)\bigr\rrvert <c_1
e^{-c_2 n}.
\end{eqnarray*}
\end{theorem}

Before we can prove Theorem \ref{theoremconvergence}, we need to
present some auxiliary definitions related to Stieltjes functions. In the
\hyperref[appendix1]{Appendix}, we collect several relevant results which
show the connections between Stieltjes functions and Pad\'e approximations.
%de4 #&#

\begin{definition}\label{defstieldef}
A \emph{Stieltjes function} is defined by the Stieltjes-integral representation
\[
f(z):= \int_{[0,\infty)} \frac{\nu(\d u)}{1 + zu},
\]
where $\nu(\d u)$ is a positive measure on $[0,\infty)$ whose support
has infinitely many different points, and which has finite moments
\begin{eqnarray*}
m_j&:=& \int_{0}^{\infty}u^j
\nu(\d u).
\end{eqnarray*}
\end{definition}

Formally, we may also express $f$ as a \emph{Stieltjes series}, which
may converge only at $0$, and has the following form:
%
%e15 #&#
%
\begin{equation}
\label{Stieltjesseries} f(z) = \sum_{j=0}^{\infty}(-z)^j
m_j.
\end{equation}
It is easy to see that the above series converges for $\llvert  z\rrvert  <R$ if and
only if $\operatorname{supp}(\nu) \subseteq[0,1/R]$. In this case we
will call $f(z)$ a Stieltjes function (or a Stieltjes series) with the
radius of convergence $R$.

\begin{pf*}{Proof of Theorem \ref{theoremconvergence}}
Let us denote $\eta(\d v)=\llvert  v\rrvert  ^3 \mu^*(\d v)$ and define
\begin{eqnarray*}
g(z)&:=&\int
_{(0,{1}/{\hat\rho}+{1}/{\rho}]} \frac
{\eta(\d(u-{1/\hat\rho}))}{1+uz}
\end{eqnarray*}
and $f(z):=zg(z)$. Note that $g(z)$ is a Stieltjes function with the
radius of convergence $R=(1/\rho+1/\hat\rho)^{-1}$, therefore,
according to Theorem \ref{thm544} in the \hyperref[appendix1]{Appendix}, the
Pad\'e approximations $g^{[n-1/n]}(z)$ converge to $g(z)$ exponentially
fast on compact subsets of $\c\setminus(-\infty,-R]$.

Changing the variable of integration $v=u-{1/\hat\rho}$ in (\ref
{formulapsiStieltjesfunction}), we obtain
%
%e16 #&#
%
\begin{eqnarray}
\label{proofthm2eq1} \psi(z)&=&az-z f \biggl(-\frac{z}{1+{z/\hat\rho}} \biggr)
\nonumber\\[-8pt]\\[-8pt]\nonumber
&=& az+\frac{z^2}{1+{z/\hat\rho}} g \biggl(-\frac{z}{1+{z/\hat\rho}} \biggr).
\end{eqnarray}
According to Theorem \ref{thm152} in the \hyperref[appendix1]{Appendix}, the
$[n/n]$ Pad\'e approximation is invariant under rational
transformations of the variable. Therefore, if $w=-z/(1+z/\hat\rho)$
and $F(z):=f(w)$ then $F^{[n/n]}(z)=f^{[n/n]}(w)$.
Theorem \ref{thm154} in the \hyperref[appendix1]{Appendix} shows that
$f^{[n/n]}(z)=z g^{[n-1/n]}(z)$.
Using these results, formula (\ref{proofthm2eq1})
and the fact that $\psi_n(z)=\psi^{[n+1/n]}(z)$ which was established
in Theorem \ref{thmmain}, we conclude that
%
%e17 #&#
%
\begin{eqnarray}
\label{proofthm2eq2} \psi_n(z)&=&\psi^{[n+1,n]}(z) = az-z
f^{[n/n]} \biggl(-\frac
{z}{1+{z/\hat\rho}} \biggr)
\nonumber\\[-8pt]\\[-8pt]\nonumber
&=&az+\frac{z^2}{1+{z/\hat\rho}} g^{[n-1/n]} \biggl(-\frac
{z}{1+{z/\hat\rho}} \biggr).
\end{eqnarray}
As we have noted above, the functions $g^{[n-1/n]}(z)$ converge to
$g(z)$ exponentially fast on compact subsets of
$\c\setminus(-\infty, -R]$, and it is easy to see that the function
$w(z)=-z/(1+z/\hat\rho)$ maps compact subsets
of $\c\setminus\{ (-\infty, -\hat\rho] \cup[\rho, \infty)\}$
onto compact subsets of $\c\setminus(-\infty, -R]$.
This fact combined with (\ref{proofthm2eq1}) and (\ref
{proofthm2eq2}) completes the proof of Theorem \ref{theoremconvergence}.
\end{pf*}

The results of Theorem \ref{thmmain} show that the Pad\'e approximant
$\psi^{[n+1/n]}(z)$ is always a Laplace exponent of a hyperexponential
process. However, as we have discussed in the \hyperref[sec1]{Introduction} (see the
discussion on page \pageref{threePadeapproximations}), there are two
other Pad\'e approximants, $\psi^{[n/n]}(z)$ and
$\psi^{[n+2/n]}(z)$, which can qualify as Laplace exponents. While we
do not have a counterexample, we believe that in general
it is not true that for all L\'evy processes $X \in{\mathcal
{CM}}(\hat\rho, \rho)$ the functions $\psi^{[n/n]}(z)$ and $\psi
^{[n+2/n]}(z)$ are Laplace exponents of hyperexponential processes.
However, more can be said under the additional assumption
that the process has one-sided jumps, and we present these results in
the next section.

%s2.2 #&#
\subsection{Approximating L\'evy processes with one-sided jumps}

In this section, we will consider separately two cases: when the
process $X$ has (i) jumps of finite variation or (ii) jumps of infinite
variation. In the first case, it is enough to consider subordinators
with zero linear drift
(if we know how to approximate such subordinators, we can always add a
linear drift and a Gaussian component later). Thus, we assume that $X
\in{\mathcal{CM}}(+\infty,\rho)$ is a subordinator with zero linear
drift, defined by the characteristic triple $(0,0,\pi)_{h\equiv0}$,
or, equivalently, by the Laplace exponent
%
%e18 #&#
%
\begin{equation}
\label{Laplaceexponentsubordinator} \psi(z)=\ln\e \bigl[ e^{z X_1} \bigr]=\int
_0^{\infty} \bigl(e^{zx}-1\bigr) \pi(x) \,\d x.
\end{equation}
We emphasize that while our definition of the Laplace exponent of a
subordinator is consistent with (\ref{LevyKhinchine}), it
differs from the classical definition $\phi(z):=-\ln\e[ \exp(-z
X_1) ]$ (see \cite{Bertoin,Kyprianou}).
The justification for this choice comes from the need to have a
consistent notation for all L\'evy processes under consideration: later
we will be approximating subordinators with spectrally positive
processes, and the formulas would be very confusing if we have
different notations for the Laplace exponents of these two objects.

%th3 #&#
%
\begin{theorem}\label{thmsubordinators}
Assume that $X \in{\mathcal{CM}}(+\infty,\rho)$ is a subordinator
defined by the characteristic triple $(0,0,\pi)_{h\equiv0}$.
Let $\psi(z)$ denote the Laplace exponent of~$X$, given by (\ref
{Laplaceexponentsubordinator}). Fix $k\in\{0,1,2\}$.
\begin{longlist}[(iii)]
\item[(i)]
Let $\{x_i\}_{1\le i \le n}$ and $\{w_i\}_{1\le i \le n}$ be the nodes
and the weights of the Gaussian quadrature with respect to the measure
$v^{2+k} \mu^*(\d v)$. Then
%
%e19 #&#
%
\begin{equation}
\label{phinkn} \psi^{[n+k/n]}(z)=\sum_{j=1}^k
\psi^{(j)}(0) \frac{z^j}{j!} + z^{k+1} \sum
_{i=1}^n \frac{w_i}{1-z x_i}.
\end{equation}
\item[(ii)]
The function $\psi^{[n+k/n]}(z)$ is the Laplace exponent of a
hyperexponential process $X^{(n)}$.
The process $X^{(n)}$ has a L\'evy measure with density,
%
%e20 #&#
%
\begin{equation}
\label{densitysubordinators} \pi_n(x):=\ind_{\{x>0\}}\sum
_{i=1}^n w_i x_i^{-2-k}
e^{-{x/x_i}},
\end{equation}
and is defined by the characteristic triple
%
%e21 #&#
%
\begin{eqnarray}
\label{chartriplesubordinator}
\cases{ \displaystyle (0,0,\pi_n)_{h\equiv0}, &\quad if $k=0$,
\vspace*{6pt}\cr
\displaystyle \Biggl(\psi'(0)-\sum_{i=1}^{n}
w_i/x_i,0,\pi_n\Biggr)_{h\equiv0}, &
\quad if $k=1$,
\vspace*{6pt} \cr
\displaystyle\Biggl(\psi'(0),\psi''(0)-2
\sum_{i=1}^{n} w_i/x_i,
\pi_n\Biggr)_{h\equiv x}, &\quad if $k=2$.}
\end{eqnarray}
The process $X^{(n)}$ is a subordinator if $k=0$ or $k=1$ (with zero
linear drift in the former case and positive linear drift in the latter
case), and $X^{(n)}$ is a spectrally positive process with a nonzero
Gaussian component if $k=2$.
\item[(iii)]
The functions $\psi^{[n+k/n]}(z)$ converge to $\psi(z)$ exponentially
fast on compact subsets of $\c\setminus[\rho,\infty)$.
\end{longlist}
\end{theorem}

Before proving Theorem \ref{thmsubordinators}, we need to establish
the following auxiliary result.

%le2 #&#
%
\begin{lemma}\label{propmomentminussumpositive}
Assume that $\nu(\d x)$ is a finite positive measure on $(0,R]$. Let
$\{x_i\}_{1\le i \le n}$ and $\{w_i\}_{1\le i \le n}$ be the nodes and the
weights of the Gaussian quadrature with respect to the measure $x\nu
(\d x)$ on $(0,R]$. Then
\[
\sum
_{i=1}^n {w_i/x_i}<
\int_{(0,R]} \nu(\d x).
\]
\end{lemma}

\begin{pf}
Consider two Stieltjes functions
\[
f(z):=\int_{(0,R]} \frac{\nu(\d x)}{1+xz},\qquad g(z):=\int
_{(0,R]} \frac{x\nu(\d x)}{1+xz}.
\]
It\vspace*{1pt} is easy to check that $f(z)=f(0)-zg(z)$. From Theorems \ref
{PadeandGauss} and \ref{thm154}
in the \hyperref[appendix1]{Appendix}, we find that
$f^{[n/n]}(z)=f(0)-zg^{[n-1/n]}(z)$ and $g^{[n-1/n]}(z)=\sum_{1\le i
\le n} w_i/(1+x_i z)$. Therefore,
%
%e22 #&#
%
\begin{equation}
\label{propproofeq1} \lim
_{z\to+\infty} f^{[n/n]}(z)=f(0)-\sum
_{i=1}^n
{w_i/x_i}.
\end{equation}
Consider the function $F(z):=(f(0)/f(z)-1)/z$. Note that $F(z) \to
-f'(0)/\break f(0)$ as $z\rightarrow0$, and that $F(z)$ is analytic in some
neighborhood of zero.
From Theorems \ref{thm153} and \ref{thm154} in the \hyperref[appendix1]{Appendix}, we obtain
\[
F^{[n-1/n]}(z)=\frac{1}{z} \biggl(\frac{m_0}{f^{[n/n]}(z)}-1 \biggr),
\]
which can be rewritten as
%
%e23 #&#
%
\begin{equation}
\label{fnnFn1n} f^{[n/n]}(z)={m_0/\bigl(1+zF^{[n-1/n]}(z)
\bigr)}.
\end{equation}
Theorem 1.3 in \cite{KJCB} tells us that $F(z)$ is also a Stieltjes
function, and since it is analytic in a neighborhood of zero, it has a
positive radius of convergence (and therefore, finite moments). Theorem
\ref{PadeandGauss} in the \hyperref[appendix1]{Appendix} implies that $\lim_{z\to+\infty} zF^{[n-1/n]}(z)$ is finite and positive. This fact
combined with
(\ref{fnnFn1n}) shows that $\lim_{z\to+\infty} f^{[n/n]}(z)$ is
strictly positive, and applying (\ref{propproofeq1}) we obtain the
statement of the lemma.
\end{pf}

\begin{pf*}{Proof of Theorem \ref{thmsubordinators}}
First, we note that since the process $X$ has jumps of finite
variation, Lemma \ref{lemmamustar} ensures that
$v^2 \mu^*(\d v)$ is a finite measure.
Formulas (\ref{defpi}) and (\ref{Laplaceexponentsubordinator})
give us
%
%e24 #&#
%
\begin{equation}
\label{formulaphi} \psi(z)=z\int
_{(0, {1}/{\rho}]}\frac{v^2 \mu^*(\d v)}{1-vz}.
\end{equation}
We will prove the case $k=2$, as the other two cases can be treated in
the same way.
We start with the identity (\ref{formulaphi}) and rewrite it in the
equivalent form
\begin{eqnarray*}
\psi(z) &=&z\int
_{(0, {1}/{\rho}]} v^2 \mu^*(\d v)+z^2\int
_{(0, {1}/{\rho}]}
v^3 \mu^*(\d v)+z^3\int
_{(0, {1}/{\rho}]}\frac{v^4 \mu^*(\d v)}{1-vz}
\\
&=& \psi'(0)z+\psi''(0)
\frac{z^2}{2}+z^3\int
_{(0, {1}/{\rho
}]}\frac{v^4 \mu^*(\d v)}{1-vz}.
\end{eqnarray*}

The result of item (i) follows from the above expression and Theorems
\ref{PadeandGauss} and \ref{thm154} in the \hyperref[appendix1]{Appendix}.

Let us prove (ii). We use Lemma \ref{propmomentminussumpositive},
from which it follows that
\[
\frac{1}{2}\psi''(0)-\sum
_{i=1}^{n} w_i/x_i=\int
 _{(0,
{1}/{\rho}]}
v^3 \mu^*(\d v)-\sum_{i=1}^{n}
w_i/x_i>0,
\]
and thus\vspace*{1pt} the coefficient of the Gaussian component is positive. Using
(\ref{LevyKhinchine}), we compute the Laplace exponent of the process
$X^{(n)}$ corresponding to the characteristic triple $(\psi'(0),\psi
''(0)-2\sum_{i=1}^{n} w_i/x_i,\pi_n)_{h\equiv x}$:
\begin{eqnarray*}
\psi_{X^{(n)}}(z)&=& \Biggl(\psi''(0)-\sum
_{i=1}^n \frac
{w_i}{x_i} \Biggr)
\frac{z^2}{2}+\psi'(0)z+z^2 \sum
 _{i=1}^n
\frac{w_i}{x_i(1-zx_i)}
\\
&=&\psi'(0)z+\psi''(0)
\frac{z^2}{2}+z^3 \sum
_{i=1}^n
\frac
{w_i}{1-zx_i}=\psi^{[n+2/n]}(z),
\end{eqnarray*}
which proves (ii). Item (iii) follows from (\ref{formulaphi}) and
Theorem \ref{thm544} in the \hyperref[appendix1]{Appendix}.
\end{pf*}

Now we consider the second class of processes with one-sided jumps:
spectrally positive L\'evy processes with jumps of infinite variation.
Again, without loss of generality we assume that there is no Gaussian component.
Our results are presented in the following theorem (the proof is
omitted, as it is identical to the proof of Theorem \ref{thmsubordinators}).

%th4 #&#
%
\begin{theorem}\label{thmspectrallypositive}
Assume that $X \in{\mathcal{CM}}(+\infty,\rho)$ is a spectrally
positive process having jumps of infinite variation and defined by the
characteristic triple $(a,0,\pi)_{h\equiv x}$. Let $\psi(z)$ be its
Laplace exponent defined by (\ref{LevyKhinchine}).
Fix $k\in\{1,2\}$.
\begin{longlist}[(iii)]
\item[(i)]
Let $\{x_i\}_{1\le i \le n}$ and $\{w_i\}_{1\le i \le n}$ be the nodes
and the weights of the Gaussian quadrature with respect to the measure
$v^{2+k} \mu^*(\d v)$. Then
%
%e25 #&#
%
\begin{equation}
\label{phinkn2} \psi^{[n+k/n]}(z)=\sum_{j=1}^k
\psi^{(j)}(0) \frac{z^j}{j!} + z^{k+1} \sum
_{i=1}^n \frac{w_i}{1-z x_i}.
\end{equation}
\item[(ii)]
The function $\psi^{[n+k/n]}(z)$ is the Laplace exponent of a
hyperexponential process $X^{(n)}$.
The process $X^{(n)}$ has a L\'evy measure with density,
\[
\pi(x):=\ind_{\{x>0\}}\sum_{i=1}^n
w_i x_i^{-2-k} e^{-{x/x_i}},
\]
and is defined by the characteristic triple
%
%e26 #&#
%
\begin{eqnarray}
\label{chartriplesubordinator} \cases{ \bigl(\psi'(0),0,\pi\bigr)_{h\equiv x},
&\quad if $k=1$,
\vspace*{6pt}\cr
\displaystyle\Biggl(\psi'(0),\psi''(0)-2
\sum_{i=1}^{n} w_i/x_i,
\pi\Biggr)_{h\equiv x}, &\quad if $k=2$.}
\end{eqnarray}
\item[(iii)]
The functions $\psi^{[n+k/n]}(z)$ converge to $\psi(z)$ exponentially
fast on compact subsets of $\c\setminus[\rho,\infty)$.
\end{longlist}
\end{theorem}

\begin{rem*}
Let us explain why we have three different approximations
in the case of subordinators and only two approximations in the case of
spectrally positive processes. For a spectrally positive process with
jumps of infinite variation, the measure $v^2 \mu^*(\d v)$ is not
finite (see Lemma \ref{lemmamustar}), thus we cannot define Gaussian
quadrature with respect to this measure and our method of proving that
$\psi^{[n/n]}(z)$ is a Laplace exponent (in Theorem \ref
{thmsubordinators}) will not work. While we do not have a
counterexample, we believe that it is not true that for any spectrally
positive process $X$ with completely monotone
jumps [and Laplace exponent $\psi(z)$] the function $\psi^{[n/n]}(z)$
is a Laplace exponent of a hyperexponential process.
\end{rem*}

%s3 #&#
\section{Explicit examples and extensions of the algorithm}\label{sectionexplicitexamples}

In this section, we pursue three goals. First, we will show how the
results of Theorems \ref{thmsubordinators} and
\ref{thmspectrallypositive} can lead to explicit formulas in the
case of Gamma subordinators and
one-sided tempered stable processes. Then we use these results to
construct explicit hyperexponential approximations
to VG, CGMY and NIG processes. Finally, we discuss several extensions
of the approximation technique described in the previous section.

The \textit{Jacobi polynomials} $P_n^{(\alpha,\beta)}(x)$ will play an
important role in this section. They are defined as follows:
%
%e27 #&#
%
\begin{equation}
\label{defJacobi} P_n^{(\alpha,\beta)}(x):=\sum
_{j=0}^n
\pmatrix{\alpha+n
\cr
n-j} \pmatrix{\alpha+\beta+n+j
\cr
j} \biggl(
\frac{x-1}2 \biggr)^j.
\end{equation}
When $\alpha>-1$ and $\beta>-1$, these polynomials satisfy the
orthogonality condition
%
%e28 #&#
%
\begin{eqnarray}
\label{Jacobiorthogonality}
&& \int_{-1}^1 P_n^{(\alpha,\beta)}(x) P_m^{(\alpha,\beta)}(x)
(1-x)^{\alpha}(1+x)^{\beta}\,\d x
\nonumber\\[-8pt]\\[-8pt]\nonumber
&&\qquad =\frac{2^{\alpha+\beta
+1}}{2n+\alpha+\beta+1} \frac{\Gamma(n+\alpha+1)}{\Gamma(n+\beta+1)}{\Gamma(n+\alpha +\beta+1)n!}
\delta_{n,m}.
\end{eqnarray}
See Section~8.96 in \cite{Jeffrey2007} for other results related to
Jacobi polynomials.

%\begin{eqnarray}
%P_n^{(\alpha,\beta)}(x)=\frac{(\alpha+1)_n}{n!} \sum _{j=0}^n
%\frac{(-n)_j (n+\alpha+\beta+1)_j}{(\alpha+1)_j j!} \left(
%\frac{1-x}2 \right)^j.
%\end{eqnarray}

%s3.1 #&#
\subsection{Example 1: Gamma subordinator}\label{subsectionGamma}

Consider a Gamma process $X$ with both mean rate and the variance rate
equal to one. In other words, $X$ is a subordinator with
zero linear drift, which has L\'evy density $\pi(x)=x^{-1}\exp(-x)$
for $x>0$, and Laplace exponent $\psi(z)=-\ln(1-z)$ [recall that we
are using (\ref{Laplaceexponentsubordinator}) as the definition of
the Laplace exponent of a subordinator]. The random variable
$X_t$ has Gamma distribution
\[
\p(X_t \in\d x)=\frac{1}{\Gamma(t)} x^{t-1}
e^{-x} \,\d x,\qquad x>0.
\]
The following proposition gives explicit results for the approximations
to $X$, described in Theorem \ref{thmsubordinators}.

%pr1 #&#
%
\begin{proposition}\label{propGamma}
Let $X$ be a Gamma process defined by the Laplace exponent $\psi
(z)=-\ln(1-z)$. Fix $k\in\{0,1,2\}$.
\begin{longlist}[(iii)]
\item[(i)]
The denominators of the Pad\'e approximants $\psi
^{[n+k/k]}(z)=p_{n,k}(z)/\break q_{n,k}(z)$ are given by
%
%e29 #&#
%
\begin{equation}
\label{qnkgamma} q_{n,k}(z)=z^{n} P_n^{(0,k)}
({2/z}-1 ).
\end{equation}
In the case $k=0$, the numerators are also given by an explicit formula
%
%e30 #&#
%
\begin{equation}
\label{pnkgamma} p_{n,0}(z)=2 \sum
_{j=0}^n
\pmatrix{n
\cr
j}^2 [ H_{n-j}-H_{j} ]
(1-z)^j,
\end{equation}
where $H_0:=0$ and $H_j:=1+1/2+\cdots+1/j$ for $j\ge1$.
\item[(ii)]
The nodes of the Gaussian quadrature described in Theorem \ref
{thmsubordinators} are given by $x_i=(y_i+1)/2$, where $y_i \in
(-1,1)$ are the roots of
the Jacobi polynomials $P_n^{(0,k)}(y)$.
\end{longlist}
\end{proposition}

\begin{pf}
We check that
\[
-\ln(1-z)=z \int_0^1 \frac{\d v}{1-zv},
\]
and comparing the above result with formula (\ref{formulaphi}) we
identity $v^2 \mu^*(\d v)=\d v$, which is just the Lebesgue measure on $(0,1)$.
The orthogonal polynomials with respect to the measure $\ind_{\{0<v<1\}
}v^k \,\d v$ are given by the shifted Jacobi polynomials
$P_n^{(0,k)}(2z-1)$. Formula (\ref{qnkgamma}) follows from this fact
and Theorems \ref{thmsubordinators} and \ref{PadeandGauss}.
Statement (ii) follows from the well-known fact that the nodes of the
Gaussian quadrature coincide with the roots of orthogonal polynomials
(see the \hyperref[appendix1]{Appendix}).

Using an equivalent representation for the Jacobi polynomials (see
formula (8.960.1) in \cite{Jeffrey2007})
\[
P_n^{(\alpha,\beta)}(x)=\sum
_{j=0}^n
\pmatrix{\alpha+n
\cr
n-j} \pmatrix{\beta+n
\cr
j} \biggl( \frac{x-1}2
\biggr)^{j} \biggl(\frac{x+1}2 \biggr)^{n-j},
\]
we find that
\[
q_{n,k}(z)=\sum
_{j=0}^n \pmatrix{k+n
\cr
n-j} \pmatrix{n
\cr
j} (1-z)^j.
\]
The above result and formula (5) in \cite{Weideman} give us the
explicit expression for $p_{n,0}(z)$ in (\ref{pnkgamma}).
\end{pf}

%s3.2 #&#
\subsection{Example 2: Tempered stable subordinator/spectrally positive process}\label{subsectiontemperedstable}

Consider a L\'evy process $X$ defined by the Laplace exponent
%
%e31 #&#
%
\begin{equation}
\label{Laplaceexponenttemperedstable} \psi(z)=\Gamma(-\alpha) \bigl((1-z)^{\alpha}-1\bigr),
\end{equation}
where $\alpha\in(0,1)\cup(1,2)$.
It is known (see formula (4.30) in \cite{Cont}) that the L\'evy
density of the process $X$ is given by
%
%e32 #&#
%
\begin{equation}
\label{pitemperedstable} \pi(x)= \ind_{\{x>0\}} x^{-1-\alpha} e^{-x}.
\end{equation}
When $\alpha\in(0,1)$ then $X$ is a subordinator with zero linear
drift, and when $\alpha\in(1,2)$ then
$X$ is a spectrally positive process with jumps of infinite variation
and zero Gaussian component.

%pr2 #&#
%
\begin{proposition}\label{proptemperedstable}
Let $X$ be a tempered stable process defined by the Laplace exponent
(\ref{Laplaceexponenttemperedstable}).
For $\alpha\in(0,1)$ ($\alpha\in(1,2)$) we fix a value of $k\in\{
0,1,2\}$ (resp., $k\in\{1,2\}$).
\begin{longlist}[(ii)]
\item[(i)]
The denominators and the numerators of the Pad\'e approximants $\psi
^{[n+k/k]}(z)=p_{n,k}(z)/q_{n,k}(z)$ are given by
%
%e33 #&#
%e34 #&#
%
\begin{eqnarray}
\label{qnktemperedstable} q_{n,k}(z)&=&z^{n} P_n^{(\alpha,k-\alpha)}
({2/z}-1 ),
\\
\label{pnktemperedstable} p_{n,k}(z)&=&\Gamma(-\alpha) \Biggl[\frac{1}{n!}
\sum
_{j=0}^{n+k} \frac{(2n+k-j)!(-n-\alpha
)_j}{j!(n+k-j)!}z^j-q_{n,k}(z)
\Biggr].
\end{eqnarray}
\item[(ii)]
The nodes of the Gaussian quadratures described in Theorems \ref
{thmsubordinators} and \ref{thmspectrallypositive} are given by
$x_i=(y_i+1)/2$, where $y_i\in(-1,1)$ are the roots of the Jacobi
polynomials $P_n^{(\alpha,k-\alpha)}(y)$.
\end{longlist}
\end{proposition}

\begin{pf}
We check that for $x>0$ and $\alpha>0$
\[
x^{-1-\alpha} e^{-x}=\frac{1}{\Gamma(1+\alpha)} \int_1^{\infty}
e^{-ux} (u-1)^{\alpha} \,\d u.
\]
The above result combined with formulas (\ref{defpi}), (\ref
{defmustar}) and (\ref{pitemperedstable}) gives us
\[
\mu(\d u) \sim\ind_{\{u>1\}}(u-1)^{\alpha} \,\d u\quad\mbox{and}\quad
v^2 \mu^*(\d v) \sim\ind_{\{0<v<1\}}v^{-\alpha}
(1-v)^{\alpha} \d v,
\]
where the symbol ``$\sim$'' means ``equal, up to a multiplicative
constant.'' This shows that
the orthogonal polynomials with\vspace*{1pt} respect to the measure $v^{2+k} \mu
^*(\d v)$ are given by the shifted Jacobi polynomials
$P_n^{(\alpha,k-\alpha)}(2z-1)$. Formula (\ref{qnktemperedstable})
follows from this fact and Theorems \ref{thmsubordinators} and \ref
{PadeandGauss}, and statement (ii) follows from the fact that the
nodes of the Gaussian quadrature coincide with the roots of orthogonal
polynomials (see the \hyperref[appendix1]{Appendix}). Formula (\ref
{pnktemperedstable}) follows from the last equation in \cite
{Iserles}, and the following fact: \textit{if} $m\ge n\ge1$ \textit{and}
$p(z)/q(z)$ \textit{is the} $[m/n]$ \textit{Pad\'e approximant to} $f(z)$, \textit{then}
$a(p(z)-q(z))/q(z)$ \textit{is the} $[m/n]$ approximant \textit{to} $a(f(z)-1)$. The
above fact is easy to deduce from Definition \ref{defpade}.
\end{pf}

%s3.3 #&#
\subsection{Approximating VG, CGMY and NIG processes}\label{subsectionapproximatingVGCGMYNIG}

The results of Propositions \ref{propGamma} and \ref{proptemperedstable} can be used to construct explicit approximations
to VG, NIG and CGMY processes. There are two methods for doing this:
(i) we can construct the process with two-sided jumps as a difference
of processes with only positive jumps or (ii) we can express the
process as a Brownian motion with drift, time-changed by a subordinator.

Let us describe the first approach using the example of the VG process
$X$ (see \cite{Madan}).
We will denote by $\Gamma(t;\mu,\nu)$ the Gamma process with mean
rate $\mu$ and variance rate $\nu$, defined by the Laplace exponent
\[
\psi_{\Gamma}(z)=-\frac{\mu^2}{\nu} \ln \biggl(1-\frac{\nu}{\mu
}z
\biggr).
\]
The Variance Gamma process is defined as the Brownian motion with drift
$\sigma W_t + \theta t$ subordinated by an independent Gamma process
$Y_t$ with mean rate one and variance rate $\nu$. The Laplace exponent
of $X$ is given by
%
%e35 #&#
%
\begin{equation}
\label{psiVG1} \psi_X(z)=\psi_{\Gamma} \biggl(\theta z+
\frac{\sigma
^2}{2}z^2 \biggr)=-\frac{1}{\nu} \ln \biggl(1-\nu
\theta z - \nu \frac{\sigma^2}{2}z^2 \biggr).
\end{equation}
Define $\mu_{p}=\frac{1}{2} \sqrt{\theta^2+2\sigma^2/\nu}+\theta
/2$ and $\mu_{n}=\mu_p-\theta$. The identity
%
%e36 #&#
%
\begin{equation}
\label{psiVG2} -\frac{1}{\nu} \ln \biggl(1-\nu\theta z - \nu
\frac{\sigma
^2}{2}z^2 \biggr)=-\frac{1}{\nu} \ln (1-
\mu_p \nu z )-\frac{1}{\nu} \ln (1+ \mu_n \nu z )
\end{equation}
allows us to write $X$ as the difference of two independent Gamma subordinators
%
%e37 #&#
%
\begin{equation}
\label{XVG} X_t=\Gamma\bigl(t;\mu_p,
\mu_p^2 \nu\bigr)-\Gamma\bigl(t;\mu_n,
\mu_n^2 \nu\bigr).
\end{equation}
In order to approximate the VG process $X$ by a hyperexponential
process, we use Proposition \ref{propGamma} and approximate each
Gamma process in (\ref{XVG}) by a hyperexponential subordinator
[equivalently, we approximate each logarithm in (\ref{psiVG2}) by a
rational function].

The same procedure works for CGMY processes. They are defined by the
Laplace exponent
%
%e38 #&#
%
\begin{equation}
\label{defpsiCGMY} \psi_X(z)= C \Gamma(-Y) \bigl[ (M- z
)^{Y}-M^Y+ (G+ z )^{Y}- G^Y
\bigr],
\end{equation}
where $Y \in(0,1) \cup(1,2)$ and all remaining parameters are
positive. We see that $X$ can be obtained as a linear drift plus a
difference of two scaled tempered stable processes with only positive
jumps. Proposition \ref{proptemperedstable} gives us an explicit
approximation to the one-sided processes [equivalently, explicit
approximations to each
power function in (\ref{defpsiCGMY})], and as a result we obtain
explicit hyperexponential approximations to general two-sided CGMY processes.

The second procedure for obtaining explicit approximations uses the
representation of the process as a Brownian motion with drift,
time-changed by a subordinator $Y$. The main idea is that we
approximate the subordinator $Y$ by a \mbox{hyperexponential} subordinator
$\tilde Y$, which we then use as a time-change process, instead of $Y$.
The following proposition ensures that the resulting approximation is
also hyperexponential.

%pr3 #&#
%
\begin{proposition}\label{proptimechange}
Assume that $\tilde Y$ is a hyperexponential subordinator and $W$ is an
independent Brownian motion. Then for all $\sigma>0$ and $a\in\r$ the
process $Z_t:=\sigma W_{\tilde Y_t}+a \tilde Y_t$ is also hyperexponential.
\end{proposition}

\begin{pf}
Denote the Laplace exponent of $\tilde Y$ as $\psi_{\tilde Y}(z)$.
Since $\tilde Y$ is hyperexponential, $\psi_{\tilde Y}(z)$ is a
rational function. It is well known that the Laplace exponent of the
subordinated process $Z$ is given by $\psi_Z(z)=\psi_{\tilde
Y}(\sigma^2 z^2/2+a z)$, therefore, it is also a rational function.
Proposition 2.1 in \cite{Jeannin} tells us that the process $Z$ has
completely monotone L\'evy density. This fact and rationality of $\psi
_Z$ prove that $Z$ is hyperexponential.
\end{pf}

As we have discussed above, the VG process can be obtained as a
Brownian motion with drift, time-changed by a Gamma process $Y$.
Proposition \ref{propGamma} gives us explicit hyperexponential
approximations to the Gamma process $Y$, therefore, from
Proposition \ref{proptimechange} we obtain explicit hyperexponential
approximations to the original VG process.

The same ideas can be applied to the NIG process (see Section~4.4.3 in
\cite{Cont}), which is defined as a Brownian motion with drift
time-changed by an inverse Gaussian subordinator $Y$, defined by
Laplace exponent $\psi_Y(z)=(1-\sqrt{1-\kappa z})/\kappa$.
Proposition \ref{proptemperedstable} gives us explicit
hyperexponential approximations to $Y$ and, therefore, we obtain
explicit hyperexponential approximations to the NIG process itself.

The approximations described above have a number of desirable features.
They are quite explicit, and the nodes of Gaussian quadratures which
are needed to compute the characteristic triples of the approximating
processes are expressed in terms of the roots of Jacobi polynomials
(for which there exist extensive tables, and which can also be computed
very easily by numerical means). The first method, based on decomposing
the process into a difference of one-sided processes, is also quite
flexible: we are free to choose the degree of the Pad\'e approximation
for each one-sided process independent of another. This may be helpful
in applications, such as when pricing down-and-out barrier options: we
may want to approximate negative jumps more accurately than positive jumps.
\label{pagesuboptimal} However, we would like to emphasize that these
approximations are not optimal, in the sense of property (iii) in Theorem
\ref{thmmain}: the general method for approximating two-sided L\'evy
processes gives a hyperexponential process with a Laplace exponent of
smaller degree (the degree of a rational function is defined as the
maximum of the degree of the numerator and denominator), which matches
more moments of the original process. In
Section~\ref{sectionnumerics}, we compare the numerical efficiency of
these two methods.

%s3.4 #&#
\subsection{Extensions of the approximation algorithm}\label{subsectionextensions}

There are two ways in which Theorems \ref{thmmain}, \ref
{thmsubordinators} and \ref{thmspectrallypositive} can be generalized.
First, there is an almost trivial (but potentially useful)
generalization, in that instead of considering the Pad\'e approximation
at $0$, we can consider the Pad\'e approximation centered at another
point $a \in(-\hat\rho, \rho)$. Then the statements of Theorems
\ref{thmmain}, \ref{theoremconvergence}, \ref{thmsubordinators}
and \ref{thmspectrallypositive} would still be true, provided that
we replace the Pad\'e approximation
$\psi^{[n+k/n]}(z)$ (centered\vspace*{1pt} at $0$) by $\psi^{[n+k/n]}(z)-\psi
^{[n+k/n]}(a)$ (centered at $a$). This fact can be easily established
using the Esscher transform, which
maps a L\'evy process $X \in{\mathcal{CM}}(\hat\rho, \rho)$
defined by Laplace exponent $\psi(z)$ into\vspace*{2pt} a process
$\tilde X \in{\mathcal{CM}}(\hat\rho+a, \rho-a)$, defined by
Laplace exponent $\tilde\psi(z)=\psi(a+z)-\psi(a)$.

The second generalization is that instead of a Pad\'e approximation one
can use a general rational interpolation, which can informally be
defined as a \textit{multi-point Pad\'e approximation}; see \cite
{donoghue1974}. The following algorithm describes how to approximate
the Laplace exponent $\psi(z)$ of a L\'evy process $X \in{\mathcal
{CM}}(\hat\rho, \rho)$.\vspace*{6pt}

\textit{A general approximation algorithm}:
\begin{longlist}[(iii)]
\item[(i)] For $k\ge1$, choose $k$ distinct points $\{z_i\}_{1\le i
\le k}$ inside the interval $(-\hat\rho, \rho)$.
\item[(ii)] Choose nonnegative integers $\{\beta_i\}_{1\le i \le k}$,
such that $k+\sum_{i=1}^k \beta_i=2n+1$ for some integer $n$.
\item[(iii)] We want to find a rational function $\tilde\psi(z)=zP(z)/Q(z)$
with $\operatorname{deg}(P)\le n$ and $\operatorname{deg}(Q) \le n$ and
which satisfies
%
%e39 #&#
%
\begin{equation}
\label{Cauchyinterpolationconditions}
\qquad\frac{\d^j}{\d z^j} \bigl({\tilde\psi(z)/z}\bigr) \bigl\llvert
_{z=z_i} = \frac{\d
^j}{\d z^j} \bigl( {\psi(z)/z} \bigr) \bigr\rrvert
_{z=z_i},\qquad 1\le i \le k, 0\le j \le\beta_i.
\end{equation}
\end{longlist}

%th5 #&#
%
\begin{theorem}\label{thmgeneralapproximation}
Assume that $\psi(z)$ is the Laplace exponent of a L\'evy process $X
\in{\mathcal{CM}}(\hat\rho,\rho)$. There exists a unique rational
function $\tilde\psi(z)$ which satisfies the conditions of item \textup{(iii)}.
Moreover, $\tilde\psi(z)$ is the Laplace exponent of a
hyperexponential process $\tilde X \in{\mathcal{CM}}(\hat\rho,\rho)$.
\end{theorem}

\begin{pf}
We recall that $f(z)$ is called a \textit{Pick function} if $f(z)$ is
analytic in the upper-half plane ${\mathbb H}=\{z \in\c\dvtx  \im
(z)>0\}$ and satisfies $f({\mathbb H}) \subseteq{\mathbb H}$. There
exists a bijection between L\'evy processes with completely monotone
jumps and Pick functions: $X\in{\mathcal{CM}}(\hat\rho, \rho)$ if
and only if $\psi(z)/z$ is a Pick function analytic in
$\c\setminus\{(-\infty, -\hat\rho] \cup[\rho,\infty)\}$ (see
\cite{Kwasnicki}, Theorem 5.1, or \cite{Rogers}).
The result of Theorem \ref{thmgeneralapproximation} now follows
easily from the fact that
$\psi(z)/z$ is a Pick function and \cite{donoghue1974}, Theorem 4,
which guarantees the existence of a rational functions $\tilde\psi
(z)$ satisfying conditions (\ref{Cauchyinterpolationconditions}) and
states that $\tilde\psi(z)/z$ is a also Pick function, which is
analytic in $\c\setminus\{(-\infty, -\hat\rho] \cup[\rho,\infty
)\}$.
\end{pf}

%f1 #&#
%
\begin{figure}[b]

\includegraphics{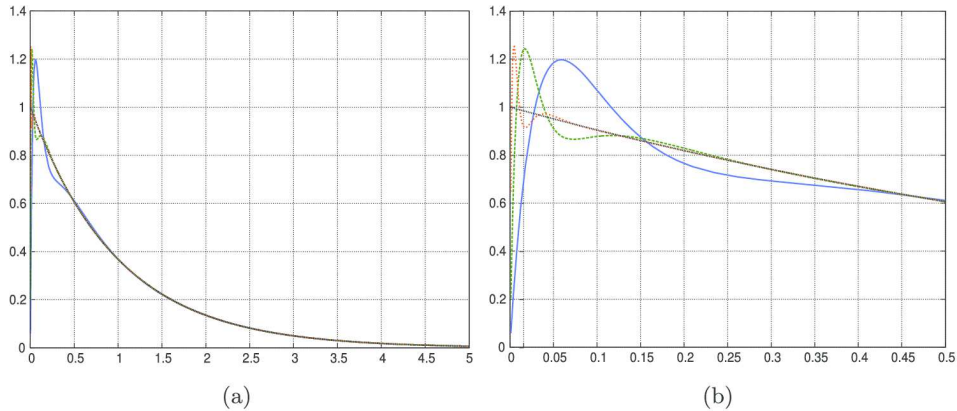}

\caption{The graph of $x \pi(x)$ (black curve) and $x \pi
^{[n/n]}(x)$, where $\pi(x)=x^{-1}\exp(-x) {\mathbf1}_{\{x>0\}}$ is
the L\'evy density of the Gamma subordinator, and $\pi^{[n/n]}(x)$ is
the L\'evy density corresponding to $\psi^{[n/n]}(z)$ Pad\'e
approximation, given by formula (\protect\ref{densitysubordinators}).
Blue, green and red curves correspond to $n
\in\{5,10, 20\}$. \textup{(b)}~depicts the magnification by a factor of
ten of the region near the origin in~\textup{(a)}.}\label{figLevydensity}
\end{figure}

%s4 #&#
\section{Numerical results}\label{sectionnumerics}

In this section we present a number of numerical experiments, which
demonstrate the efficiency of our approximations.
As a first example, we consider the Gamma process $X$ defined by the
Laplace exponent $\psi(z)=-\ln(1-z)$. We compute the L\'evy density
$\pi_n(x)$ corresponding to the approximation $\psi^{[n/n]}(z)$,
which is given explicitly in Proposition \ref{propGamma}. The L\'evy
density of the Gamma process is given by $\pi(x)=\exp(-x)/x$, thus in
order to avoid the singularity at $x=0$ we compare the graphs
of $x\pi(x)\equiv\exp(-x)$ and $x\pi_n(x)$. The results are
presented on Figure~\ref{figLevydensity}. We see that even with a
small value of $n=5$ the tail of $\pi_n(x)$ matches the tail of $\pi
(x)$ very well, and as $n$ increases the approximation converges very
rapidly (as long as $x$ is not too close to zero).

Next, we compare the cumulative distribution function (CDF) of $X_t$
for the same Gamma process $X$ and its approximations $X^{(n,k)}$,
which are defined by the Laplace exponents $\psi^{[n+k/n]}(z)$, $k \in
\{0,1,2\}$, see Proposition \ref{propGamma}.
We compute the CDF for two values of $t \in\{1,2\}$. The CDF of $X_t$
for the Gamma process is known explicitly:
\[
\p(X_1\le x)=1-e^{-x}\quad\mbox{and}\quad
\p(X_2 \le x)=1-(x+1)e^{-x}. %
\]
We define the numbers $r_i$ as the coefficients in the
asymptotic expansion
\[
\psi^{[n+k/n]}(z)=r_2 z^2 + r_1 z +
r_0 + O(1/z),\qquad z\to\infty, %
\]
and define
\begin{eqnarray*}
\phi_{n,k}(z):= \cases{ e^{t \psi^{[n/n]}(z)}-e^{\tr_0}, &\quad if
$k=0$,
\cr
e^{t \psi^{[n+1/n]}(z)}-e^{\tr_0+\tr_1 z}, &\quad if $k=1$,
\cr
e^{t \psi^{[n+2/n]}(z)}, &\quad if $k=2$.}
\end{eqnarray*}
The CDF of the approximating process is computed by the Fourier inversion
%
%e40 #&#
%
\begin{eqnarray}
\label{formulaCDF} \p\bigl(X^{(n,k)}_t \le x\bigr)&=&1-
\frac{e^{-c x}}{\pi} \re \biggl[ \int_0^{\infty}
\phi_{n,k}(c+\i u) e^{-\i u x} \frac{\d u}{c+\i u} \biggr]
\nonumber\\[-8pt]\\[-8pt]\nonumber
&&{}- e^{\tr_0} \ind_{\{x \le r_1 t\}} \ind_{\{k=1\}},
\end{eqnarray}
where $x>0$ and $c \in(0,1)$.

Let us explain the intuition behind the formula corresponding to $k=0$,
the other cases can be treated similarly. According to Theorem \ref
{thmsubordinators}, the process $X^{(n,0)}$ is a compound Poisson
hyperexponential process with intensity $-r_0$, thus its distribution
has an atom at zero: $\p(X^{(n,0)}_t=0)=\exp(\tr_0)$. If we subtract
the atom at zero, we obtain an absolutely continuous
positive measure
%
%e41 #&#
%
\begin{equation}
\label{definenut} \nu_t(\d x):=\p\bigl(X^{(n,0)}_t
\in\d x\bigr)-e^{\tr_0} \delta_0(\d x)
\end{equation}
which has Fourier transform
\[
\int_{\r} e^{\i tz }\nu_t (\d
x)=e^{t\psi^{[n/n]}(z)}-e^{\tr_0}=\phi_{n,0}(z).
\]
Since $\nu_t(\d x)$ is absolutely continuous with total mass $1-\exp
(\tr_0)$,
we can find the CDF corresponding to this measure by the inverse
Fourier transform
%
%e42 #&#
%
\begin{equation}
\label{formulanut} \nu_t\bigl((0,x)\bigr)=1-e^{\tr_0}-
\frac{e^{-c x}}{\pi} \re \biggl[ \int_0^{\infty}
\phi_{n,0}(c+\i u) e^{-\i u x} \frac{\d u}{c+\i u} \biggr].
\end{equation}
Note that the integral in (\ref{formulanut}) converges absolutely,
since $\phi_{n,0}(c+\i u)=O(1/u)$ as $u \to\infty$.
Formula (\ref{formulaCDF}) follows at once from (\ref{definenut})
and (\ref{formulanut}).

The results of our computations are presented in Table~\ref{tabCDF}.
We see that the CDF of $X^{(n,k)}_t$ does converge to
$X_t$, and the convergence seems to be faster for $t=2$ than it is for $t=1$.

%t1 #&#
%
\begin{table}
\tabcolsep=0pt
\caption{The values of $\varepsilon_{n,k}(t):=\max_{x\ge0} \llvert  \p(X_t\le
x) - \p(X_t^{(n,k)}\le x)\rrvert$, where $X$ is the Gamma process with $\psi
(z)=-\ln(1-z)$ and the process $X^{(n,k)}$ has Laplace exponent $\psi
^{[n+k/n]}$}\label{tabCDF}
\begin{tabular*}{\tablewidth}{@{\extracolsep{\fill}}l ccc@{\qquad}c ccc@{}}
\hline
$\bolds{\varepsilon_{n,k}(1)}$ & $\bolds{k=0}$ & $\bolds{k=1}$ & $\bolds{k=2}$  &   $\bolds{\varepsilon_{n,k}(2)}$ & $\bolds{k=0}$ & $\bolds{k=1}$ & $\bolds{k=2}$ \\
\hline
$n=5$ & 1.1e--2 & 1.1e--2 & 8.8e--3                                       &   $n=5$ & 3.3e--4 & 3.2e--4 & 5.4e--4 \\
$n=10$ & 2.8e--3 & 3.4e--3 & 2.8e--3                                      &   $n=10$ & 2.6e--5 & 2.8e--5 & 5.6e--5 \\
$n=15$ & 1.3e--3 & 1.6e--3 & 1.4e--3                                      &   $n=15$ & 5.4e--6 & 6.4e--6 & 1.3e--5 \\
$n=20$ & 7.5e--4 & 9.3e--4 & 8.1e--4                                      &   $n=20$ & 1.8e--6 & 2.1e--6 & 4.6e--6 \\
\hline
\end{tabular*}
\end{table}

Our remaining examples are all related to pricing European and various
exotic options in L\'evy driven models. We will work with the following
two processes: the VG process $V$ defined by the Laplace exponent
\[
\psi(z) = \mu z -\frac{1}{\nu}\ln \biggl(1 - \frac{z}{a} \biggr) -
\frac{1}{\nu}\ln \biggl(1 + \frac{z}{\hat{a}} \biggr),
\]
and parameters
\[
(a,\hat{a},\nu) = (21.8735, 56.4414, 0.20),
\]
and the CGMY process $Z$ defined by the Laplace exponent
\[
\psi(z) = \mu z+C\Gamma(-Y) \bigl[(M - z)^Y - M^Y + (G
+ z)^Y - G^Y \bigr],
\]
and parameters
\[
(C,G,M,Y) = (1, 8.8, 14.5, 1.2).
\]
Note that $V$ is a process with jumps of infinite activity and finite
variation, whereas $Z$ has jumps of infinite variation. Both of these
processes have zero Gaussian component.
The process $V$ with the same parameters was considered in \cite
{barbench}, and later we will use their numerical results as a
benchmark for our computations.

Our approach from here on is to compare a benchmark option price (for a
variety of options) with a price calculated using one of four possible
approximations. The first approximation is based on the $[n+1/n]$ Pad\'
e approximant for the process with two-sided jumps from Theorem \ref
{thmmain}. The other three approximations are based on the algorithm
presented in Section~\ref{subsectionapproximatingVGCGMYNIG}, which
considers the process as a difference of two processes having only
positive jumps, and uses the explicit $[N+k/N]$ Pad\'e approximations
from Propositions \ref{propGamma} and \ref{proptemperedstable}.
Note that the first approximation will result in a rational function of
degree $n+1$, while the other three approximations result in a
rational function of degree $2N+k$. In instances where we calculate
multiple approximations, we set $n=2N$ in order to make a fair
comparison between different approximations. In all examples, we define
the stock price process as $S_t=S_0 \exp(X_t)$ (where $X\equiv V$ in
the VG case or $X\equiv Z$ in the CGMY case). Further, we choose the
value of the linear drift $\mu$ such that the process $S_t \exp(-rt)$ is a martingale.

%t2 #&#
%
\begin{table}[t]
\tabcolsep=0pt
\caption{The error in computing the price of the European call option
for the VG $V$-model. The benchmark price is 2.5002779303}\label{tabVGeuropean}
\begin{tabular*}{\tablewidth}{@{\extracolsep{\fill}}@{}lcccc@{}}
\hline
& \multicolumn{1}{c}{\textbf{Two-sided}} & \multicolumn{3}{c@{}}{\textbf{One-sided}}\\[-6pt]
& \multicolumn{1}{c}{\hrulefill} & \multicolumn{3}{c@{}}{\hrulefill}\\
& $\bolds{[2N+1/2N]}$ & $\bolds{[N/N]}$ & $\bolds{[N+1/N]}$ & $\bolds{[N+2/N]}$\\
\hline
$N=1$ & $-$1.58e--2 & \phantom{$-$}9.12e--2 & \phantom{$-$}7.02e--3 & $-$3.02e--2 \\
$N=2$ & \phantom{$-$}1.66e--3 & $-$6.16e--3 & \phantom{$-$}4.80e--3 & $-$7.82e--4 \\
$N=3$ & \phantom{$-$}6.20e--4 & $-$1.28e--3 & $-$4.32e--5 & \phantom{$-$}6.78e--4 \\
$N=4$ & \phantom{$-$}1.25e--4 & \phantom{$-$}1.88e--4 & $-$1.98e--4 & \phantom{$-$}9.81e--5 \\
$N=5$ & $-$7.19e--5 & \phantom{$-$}8.82e--5 & $-$2.62e--5 & $-$2.40e--5 \\
$N=7$ & \phantom{$-$}4.34e--6 & $-$8.48e--6 & \phantom{$-$}5.82e--6 & $-$1.71e--6 \\
$N=9$ & $-$7.72e--8 & \phantom{$-$}3.31e--7 & $-$6.99e--7 & \phantom{$-$}7.35e--7 \\
$N=12$ & \phantom{$-$}4.85e--7 & $-$1.81e--8 & \phantom{$-$}4.97e--8 & $-$6.10e--8 \\
$N=15$ & $-$8.56e--8 & $-$1.37e--9 & $-$3.31e--9 & \phantom{$-$}6.06e--9 \\
\hline
\end{tabular*}
\end{table}

%t3 #&#
%
\begin{table}[b]
\tabcolsep=0pt
\caption{The error in computing the price of the European call option
for the CGMY $Z$-model. The benchmark price is 11.9207826467}\label{tabCGMYeuropean}
\begin{tabular*}{\tablewidth}{@{\extracolsep{\fill}}@{}lccc@{}}
\hline
& \multicolumn{1}{c}{\textbf{Two-sided}} & \multicolumn{2}{c@{}}{\textbf{One-sided}}\\[-6pt]
& \multicolumn{1}{c}{\hrulefill} & \multicolumn{2}{c@{}}{\hrulefill}\\
& $\bolds{[2N+1/2N]}$ & $\bolds{[N+1/N]}$ & $\bolds{[N+2/N]}$ \\
\hline
$N=1$ & $-$2.75e--2 & \phantom{$-$}1.93e--2 & $-$3.72e--3 \\
$N=2$ & $-$4.86e--6 & $-$4.19e--6 & \phantom{$-$}9.5e--5\phantom{0} \\
$N=3$ & \phantom{$-$}4.80e--7 & $-$1.48e--5 & $-$2.54e--7 \\
$N=4$ & \phantom{$-$}2.9e--8\phantom{0} & \phantom{$-$}6.41e--7 & $-$1.55e--7 \\
$N=5$ & \phantom{$-$}1.14e--9 & \phantom{$-$}5.58e--9 & \phantom{$-$}6.95e--9 \\
\hline
\end{tabular*}
\end{table}

First, we compute the price of a European call option with $S_0=100$,
strike price $K=100$, maturity $T=0.25$ and interest rate $r=0.04$. All
option prices are computed using the Fourier transform approach from
\cite{Carr}. When dealing with hyperexponential processes, we have
slightly modified this approach by removing the possible atom in the
distribution of $X_t^{(n,k)}$,
in the same way as we did earlier in equation (\ref{formulaCDF}). The
benchmark prices for the original VG process $V$ and the CGMY process
$Z$ were computed multiple times, with different discretizations of the
Fourier integral, and seem to be correct to at least $\pm$1.0e--9. The
results of our computations for the approximations to VG (CGMY) model
are presented in Table~\ref{tabVGeuropean} (resp., \ref
{tabCGMYeuropean}). We see that all four approximations are doing an
excellent job, and already for $N=4$ we obtain acceptable accuracy of
around 1.0e--4. We would like to point out that the three approximations
based on explicit one-sided approximations have remarkably good
accuracy. As we have discussed on page~\pageref{pagesuboptimal},
these approximations are not optimal in the sense that one can find a
rational Laplace exponent of lower degree which matches more moments of
the original process. However, this nonoptimality does not seem to play
any role here. These three one-sided approximations are superior to the
two-sided approximation, in the sense that they have very good accuracy
\emph{and} explicit formulas.

We also note that all four approximations seem to be doing a better job
in the case of the CGMY process $Z$. We think that the likely cause is
that the process $Z$ has jumps of infinite variation and $Z_t$ has
smooth density, which is not the case for the process $V$.

Next, we compute the price of a continuously sampled arithmetic Asian
call option with fixed strike. That is, we calculate the following quantity:
%
%e43 #&#
%
\begin{equation}
\label{eqasia} C(S_0,K,T):=e^{-r T} \e \biggl[ \biggl(\int
_{0}^T S_u \,\d u - K \biggr)^+
\biggr].
\end{equation}
We set the parameters as follows:
\[
S_0 = 100,\qquad r=0.03,\qquad T=1 %
\]
and $K=90$ for the VG process and $K=110$ for the CGMY process.
In order to compute the price of the Asian option, we use the technique
pioneered for hyperexponential processes by Cai and Kou \cite
{CaiKouAsianoptions} (see also \cite{Aopt}, Section~4.2). Since we
were unable to find any results in the literature for pricing such
options for either the VG or CGMY process (other than by Monte Carlo
methods), we use our own benchmark calculated using a significantly
larger $N$.
By experimenting with different ways of discretizing the resulting
integrals in the inverse Laplace and inverse Mellin transform
(see \cite{Aopt}, Section~4.2), we arrive at a benchmark price of
11.18859 for the process $V$
and 9.95930 for the process $Z$. These benchmark prices seem to be
correct to within $\pm$1.0e--5. The results for each $N$ are compared
to the benchmark price, the errors are gathered in Table~\ref
{tabvgasia} for the process $V$ and in Table~\ref{tabcgmyasia} for
the process $Z$.

We observe again, that convergence to the benchmark price is very rapid
and that there is little difference in the rate of convergence between
the one-sided and two-sided approximations. We note that we achieve an
acceptable error of \mbox{$\pm$1.0e--4} with a rational approximation of
degree $5$. We would like to emphasize that the numbers in Table~\ref
{tabvgasia} and Table~\ref{tabcgmyasia} represent the difference
between the approximate price and the benchmark price, and the
benchmark itself is only an approximation to the exact price with
accuracy of the order of $\pm$1.0e--5. The two most important factors
influencing the accuracy of the benchmark price are (i) the accuracy of
the approximation of the target L\'evy process $X$ by a
hyperexponential process $X^{(n)}$, and (ii) the error in the
discretization of the inverse Laplace and inverse Mellin transform
needed to compute the price of Asian option (see \cite{Aopt}, Section~4.2) in the model driven by the process $X^{(n)}$. The
results presented above only measure the effect of the first of these
two factors, which explains why some numbers are of the order $\pm
$1.0e--7 whereas our benchmark is only correct to within $\pm$1.0e--5.

%t4 #&#
%
\begin{table}
\tabcolsep=0pt
\caption{The error in computing the price of the Asian option for the
VG $V$-model. The benchmark price is 11.188589
(calculated using the $[91/90]$ two-sided approximation)}\label{tabvgasia}
\begin{tabular*}{\tablewidth}{@{\extracolsep{\fill}}@{}lcccc@{}}
\hline
& \multicolumn{1}{c}{\textbf{Two-sided}} & \multicolumn{3}{c@{}}{\textbf{One-sided}}\\[-6pt]
& \multicolumn{1}{c}{\hrulefill} & \multicolumn{3}{c@{}}{\hrulefill}\\
& $\bolds{[2N+1/2N]}$ & $\bolds{[N/N]}$ & $\bolds{[N+1/N]}$ & $\bolds{[N+2/N]}$ \\
\hline
$N=1$ & $-$1.87e--3 & \phantom{$-$}1.01e--3 & $-$1.82e--3 & \phantom{$-$}9.88e--4 \\
$N=2$ & \phantom{$-$}9.49e--5 & \phantom{$-$}2.89e--4 & $-$6.33e--5 & \phantom{$-$}3.27e--5 \\
$N=3$ & \phantom{$-$}1.30e--6 & \phantom{$-$}8.85e--6 & $-$4.24e--6 & \phantom{$-$}3.99e--6\\
$N=4$ & $-$2.83e--6 & \phantom{$-$}1.07e--6 & $-$1.36e--6 & \phantom{$-$}3.16e--7 \\
$N=5$ & $-$1.11e--7 & $-$2.48e--8 & $-$5.91e--7 & $-$3.81e--7\\
\hline
\end{tabular*}
\end{table}

%t5 #&#
%
\begin{table}[b]
\tabcolsep=0pt
\caption{The error in computing the price of the Asian option for the
CGMY $Z$-model. The benchmark price is 9.959300 (calculated using the
$[91/90]$ two-sided approximation)}\label{tabcgmyasia}
\begin{tabular*}{\tablewidth}{@{\extracolsep{\fill}}@{}lccc@{}}
\hline
& \multicolumn{1}{c}{\textbf{Two-sided}} & \multicolumn{2}{c@{}}{\textbf{One-sided}}\\[-6pt]
& \multicolumn{1}{c}{\hrulefill} & \multicolumn{2}{c@{}}{\hrulefill}\\
& $\bolds{[2N+1/2N]}$ & $\bolds{[N+1/N]}$ & $\bolds{[N+2/N]}$ \\
\hline
$N=1$ & \phantom{$-$}1.88e--4 & \phantom{$-$}7.42e--4 & $-$1.19e--3 \\
$N=2$ & \phantom{$-$}4.03e--6 & \phantom{$-$}9.05e--5 & \phantom{$-$}5.39e--6 \\
$N=3$ & $-$3.58e--7 & $-$2.64e--6 & \phantom{$-$}7.93e--8\\
$N=4$ & $-$3.88e--7 & $-$1.01e--7 & $-$1.21e--7 \\
$N=5$ & $-$5.26e--7 & $-$2.47e--7 & $-$2.49e--7\\
\hline
\end{tabular*}
\end{table}

Our final example is related to pricing down-and-out barrier put
option. That is, we wish to calculate
\[
D(S_0,K,B,T):= e^{-r T} \e \bigl[ (K-S_T )^+
{\mathbf1}_{\{
S_t > B~\mathrm{for}~0\le t \le T\}} \bigr],
\]
where $B$ is the barrier level. We calculate barrier option prices for
the process $V$, for four values
$S_0 \in\{81, 91, 101, 111\}$ and with other parameters given by
$K=100$, $B=80$, $r=0.04879$ and $T=0.5$. We use the prices computed
in \cite{barbench} as the benchmark (these prices seem to be accurate
to about $\pm$1.0e--3). In order to compute the prices of down-and-out
put options for hyperexponential processes, we use the Laplace
transform inversion method by Jeannin and Pistorius \cite{Jeannin}.
In this case, we present the results only for the one-sided $[N+1/N]$
approximations. The results are presented in Table~\ref{tabvgbar}.
We see that in almost all cases the convergence is very rapid, and we
are able to match the first four digits of the benchmark price. The
convergence is somewhat slower for $S_0=81$, which is to be expected:
it is known that when pricing barrier options, the behavior of the
price near the barrier is very sensitive to the nature of the small
jumps of the underlying process (see~\cite{BoInLe}). Therefore, we may
expect that our results will not be very precise when $S_0$ is close to
$B$, since we are approximating a process with jumps of infinite
activity by a compound Poisson process with drift.

%t6 #&#
%
\begin{table}
\tabcolsep=0pt
\caption{Barrier Option prices calculated for the VG process
$V$-model. Benchmark prices were obtained from \cite{barbench},
Table~4, column 2}\label{tabvgbar}
\begin{tabular*}{\tablewidth}{@{\extracolsep{\fill}}@{}lcccc@{}}
\hline
& $\bolds{S_0=81}$ & $\bolds{S_0=91}$ & $\bolds{S_0=101}$ & $\bolds{S_0 =111}$ \\[-6pt]
& \multicolumn{1}{c}{\hrulefill} & \multicolumn{1}{c}{\hrulefill} & \multicolumn{1}{c}{\hrulefill} & \multicolumn{1}{c@{}}{\hrulefill}\\
%\hline
\textbf{Benchmark} & \textbf{3.39880} & \textbf{7.38668} & \textbf{1.40351} & \textbf{0.04280} \\
\hline
$N=2$ & 3.44551 & 7.39225 & 1.40527 & 0.04233\\
$N=4$ & 3.40209 & 7.38957 & 1.40329 & 0.04258\\
$N=6$ & 3.39910 & 7.38939 & 1.40332 & 0.04258\\
$N=8$ & 3.39856 & 7.38936 & 1.40332 & 0.04258\\
$N=10$ & 3.39853 & 7.38936 & 1.40332 & 0.04258\\
\hline
\end{tabular*}
\end{table}

Let us describe the computing environment used for our numerical
experiments. The code was written in Fortran-90, and we used a standard
2011 laptop (with an Intel Core i5-2540M CPU). All Pad\'e
approximations were computed using the most basic algorithm based on
solving the system of linear equations
(\ref{findb}), which we describe below in the \hyperref[appendix1]{Appendix}.
Since this system of linear equations is typically ill-conditioned, all
computations related to Pad\'e approximations were performed with a
high precision of 200 digits, using the MPFUN multiple precision package
\cite{Bailey95afortran90}. The computation time of the Pad\'e
approximation and its partial fraction decomposition was on the order
of 0.1 seconds. Our goal in this section was to demonstrate the
accuracy of Pad\'e-based hyperexponential approximations, therefore, we
did not try to write the most efficient code for computing the option
prices. However, our computations were reasonably fast: the computation
time for a single European (resp., Asian or barrier) option price was
around 2 (resp., 5 or 15) seconds.

%s5 #&#
\section{Concluding remarks}\label{sectionremarks}

As we have mentioned in the \hyperref[sec1]{Introduction}, there exist other methods for
approximating processes with completely monotone jumps by
hyperexponential processes. The first of these was proposed by Jeannin
and Pistorius in \cite{Jeannin}, and the second one by Crosby, Le Saux
and Mijatovi\'c in \cite{Crosby}. Our research in this field was
initially inspired by these two papers, and we would like to summarize
their methods and highlight the similarities and differences with our method.

The approach of Jeannin and Pistorius is essentially based on
minimizing the $L_2$ distance between the target L\'evy density $\pi
(x)$ and the approximating hyperexponential L\'evy density $\pi_n(x)$.
More precisely, we are looking for a hyperexponential L\'evy density
$\pi_n(x)$ of the form
(\ref{defhyperexponentialpi}) which minimizes
%
%e44 #&#
%
\begin{equation}
\label{JPapproach} \Delta_{n,\varepsilon}=\int_{\r\setminus[-\varepsilon, \varepsilon]} \bigl(\pi
(x) - \pi_n(x)\bigr)^2 \,\d x,
\end{equation}
where $\pi(x)$ is the target L\'evy density of a process with
completely monotone jumps. Note that we do need to remove an $\varepsilon
$-neighborhood of zero in the domain of integration in (\ref
{JPapproach}), because otherwise the integral may not converge.
According to the definition of $\pi_n(x)$ in (\ref
{defhyperexponentialpi}), the quantity $\Delta_{n,\varepsilon}$ can be
considered as a function of $2N+2\hat N$ parameters
$\{\alpha_i, \beta_i\dvtx  1\le i \le N\}$ and $\{\hat\alpha
_i, \hat\beta_i\dvtx  1\le i \le\hat N\}$,
and ideally one would try to find the absolute minimum of this function
in order to get the best fit of the hyperexponential L\'evy density
$\psi_n(x)$ to the target density $\pi(x)$. Since this optimal
approach would result in a complicated nonlinear minimization problem,
it is much easier to fix the parameters $\beta_i$ and
$\hat\beta_i$ [which specify the exponents of the exponential
functions in (\ref{defhyperexponentialpi})] and to minimize over the
remaining parameters $\alpha_i$ and $\hat\alpha_i$. This
simplification results in a simpler linear problem, which can be easily
solved numerically.

Next, let us summarize the main ideas behind the method of Crosby, Le
Saux and Mijatovi\'c \cite{Crosby}.
We start with the L\'evy process with completely monotone jumps and
zero Gaussian component. We use formula (\ref
{formulapsiStieltjesfunction0}), choose a parameter $A>0$ large
enough and derive the following approximation:
\begin{eqnarray*}
\psi(z)&=& a z + z^2 \int_{\r} \frac{{\operatorname{sign}}(u)}{u-z}
\frac{\mu(\d u)}{u^2}
\\
&=& a z +z^2\int_{\r\setminus[-A,A]} \frac{{\operatorname
{sign}}(u)}{u-z}
\frac{\mu(\d u)}{u^2} + z^2 \int_{[-A,A]}
\frac{{\operatorname{sign}}(u)}{u-z} \frac{\mu
(\d u)}{u^2}
\\
&=& a z + z^2\int_{\r\setminus[-A,A]} \frac{1}{1-z/u}
\frac{\mu(\d
u)}{\llvert  u\rrvert  ^3}+z^2 \int_{[-A,A]}
\frac{{\operatorname{sign}}(u)}{u-z} \frac{\mu(\d u)}{u^2}
\\
& \approx& a z + z^2\int_{\r\setminus[-A,A]}
\frac{\mu(\d
u)}{\llvert  u\rrvert  ^3}+ z^2\int_{[-A,A]}
\frac{{\operatorname{sign}}(u)}{u-z} \frac{\mu(\d u)}{u^2}=:\tilde\psi(z),
\end{eqnarray*}
where in the last step we used the fact that $\llvert  u\rrvert  >A\gg1$ and,
therefore, $1-z/u$ can be approximated by $1$. The above approximation
is the first step in the method of Crosby et al., and it gives us
the Laplace exponent of a L\'evy process $\tilde X$ with a small (but
nonzero) Gaussian coefficient
\[
\sigma^2=2 \int_{\r\setminus[-A,A]} \frac{\mu(\d u)}{\llvert  u\rrvert  ^3}.
\]
The process $\tilde X$ has L\'evy measure $\tilde\pi(x)$, given by
(\ref{defpi}) with $\mu(\d x)$ replaced by
$\mu(\d x) {\mathbf1}_{\{\llvert  x\rrvert  \le A\}}$. It is easy to see that $\tilde\pi
(x)$ is a finite measure, thus $\tilde X$ has compound Poisson jumps.
Intuitively, the effect of this first step is to replace the jumps of
$X$ (which could be of infinite activity or infinite variation) by
compound Poisson jumps and a small Gaussian component. The second step
in the method of Crosby et al. consists in discretizing the integral
\[
\int_{[-A,A]} \frac{{\operatorname{sign}}(u)}{u-z} \frac{\mu(\d
u)}{u^2}\approx\sum
\frac{{\operatorname{sign}}(x_i)}{x_i-z} \frac{w_i}{x_i^2} %
\]
via the Gauss--Legendre quadrature (a Gaussian quadrature on the
interval $[-A,A]$ with respect to the Lebesgue measure). Combining
these two steps results in a Laplace exponent of approximating
hyperexponential process.

Our method is quite similar to the approach of Crosby, Le Saux and
Mijatovi\'c.
Instead of their first approximating step, we perform a change of
variables $u=1/v$ in the integral (\ref
{formulapsiStieltjesfunction0}). This simple trick and Assumption \ref{ass1}
give us a finite domain of integration in the $v$-variable in (\ref
{formulapsiStieltjesfunction}), so that we can apply Gaussian
quadrature with respect to the measure $\llvert  v\rrvert  ^3\mu^*(\d v)$. It turns
out that this seemingly small modification has profound consequences.
First of all, we do not need to truncate the integrals and we do not
require any external parameters (such as $\varepsilon$ or $A$ in the
above two methods). Second, our approximating Laplace exponents $\psi
_n(z)$ have a simple analytic interpretation as Pad\'e approximations
of the target Laplace exponent $\psi(z)$, which allows us to borrow
tools and ideas from the well developed theory of rational
approximations and orthogonal polynomials. Third, our approximation
turns out to be optimal in the sense that the hyperexponential process
$X^{(n)}$ constructed in Theorem \ref{thmmain} matches $2n+1$ moments
of the target process $X$
[see the statement of Theorem \ref{thmmain}(iii)]. Note that this is
the best that one can hope for: according to formula (\ref{defpsin})
the process $X^{(n)}$ has $2n+1$ free parameters, thus we cannot expect
to be able to match more than $2n+1$ moments of $X$. Finally, we show
in Theorem \ref{theoremconvergence} that our approximations converge
exponentially in $n$, where $n$ is the number of terms in the L\'evy
density, and this fast convergence is confirmed by our numerical experiments.

In conclusion, we would like to discuss how our current results fit in
the context of recent developments on meromorphic processes \cite
{Kuz2010a,Kuz2010b,KuzKyPa2011}. The main motivation for introducing
meromorphic processes was the perceived lack of explicit examples of L\'
evy processes which would be useful for modeling purposes and
convenient for numerical calculations. Meromorphic processes serve this
purpose quite well: they are flexible enough to allow for jumps of
infinite activity or infinite variation, they have many parameters and
are very similar to the widely used CGMY and VG processes, and at the
same time, they are analytically tractable and enjoy an explicit
Wiener--Hopf factorization. A meromorphic process can be informally
defined as a hyperexponential process with infinitely many terms in the
L\'evy density [so that the L\'evy measure is given by (\ref
{defhyperexponentialpi}) with the finite sum replaced by infinite
series]. Hyperexponential processes can be considered as a subclass of
meromorphic processes, in the same way that rational functions can be
considered a subclass of meromorphic functions. This turns out to be a
very useful analogy, and it seems that every formula related to
hyperexponential processes has a corresponding analogue for meromorphic
processes, with the only change that the finite sums or products would
be replaced by appropriate infinite series or products.
While hyperexponential processes are much simpler objects to work with,
compared with meromorphic processes, their big disadvantage is that
they do not allow for jumps of infinite activity or infinite variation.

As an example of how our current work complements the previous
developments on meromorphic processes, consider the following
hypothetical situation. Suppose that we have data on European options
for a certain stock and we want to price barrier options on the same
underlying stock and we want to use the CGMY model to describe the
stock price dynamics. We face a problem in that numerical computation
of prices of barrier options is not so easy in the CGMY model:
algorithms based on Monte Carlo technique are not very accurate and
rather time consuming, whereas analytical methods \cite{Jeannin} are
not possible since we do not know the Wiener--Hopf factors of a CGMY
process. One way to solve this problem would be to use meromorphic
processes. We would just replace the family of CGMY processes by a very
similar family of beta-processes \cite{Kuz2010a}, and then calibrate
the parameters of a beta-process to the available data on European
options and price barrier options in the model driven by a beta-process
using the results of \cite{Kuz2010a,KuzKyPa2011}. Another way would be
to adhere to our original choice of the CGMY model: we would calibrate
the parameters of a CGMY process to the available data and then
approximate the calibrated CGMY process with hyperexponential processes
as described in this paper.
The prices of barrier options for hyperexponential processes can be
computed easily \cite{Jeannin}. It is not clear which of these two
approaches would be a better solution in practice. The first one
requires that we abandon the CGMY model and instead use meromorphic
processes, plus we have to be careful with truncating infinite products
and sums when doing numerical computations for meromorphic processes.
The second approach allows us to keep our favorite CGMY model and
simplifies the numerical computations (dealing with hyperexponential
processes is easier compared to meromorphic processes). The downside of
the second approach is that we introduce a new source of error when we
approximate a CGMY process by a hyperexponential process. Overall, we
feel that both approaches have merit and that they deserve further
investigation.

\begin{appendix}\label{appendix1}
%%%****************************************************************************************************************
%%%****************************************************************************************************************
%%%****************************************************************************************************************
%s6 #&#
\section*{Appendix: Gaussian quadrature, Pad\'e approximations and Stieltjes functions}
Consider a finite positive measure $\nu(\d x)$ on an interval $[0,a]$.
The main idea behind Gaussian quadrature is that we want to find a
measure $\tilde\nu(\d x)$, supported on $n$ points inside $[0,a]$,
which matches the first $2n-1$ moments of $\nu(\d x)$. Thus,
the weights $\{w_i\}_{1\le i \le n}$ and the nodes $\{x_i\}_{1\le i \le
n}$ of the Gaussian quadrature are uniquely defined by equations
\[
\int_{[0,a]} x^k \nu(\d x)=\sum
 _{i=1}^n
x_i^k w_i,\qquad k=0,1,\dots,2n-1.
\]
Let $\{p_n(x)\}_{n\ge0}$ be the sequence of orthogonal polynomials
with respect to the measure $\nu(\d x)$:
$\deg(p_n)=n$ and $(p_n,p_m)_{\nu}:=\int_{[0,a]} p_n(x) p_m(x) \nu
(\mathrm{d} x)= d_n \delta_{n,m}$. It is known \cite{szego}, Theorems 3.4.1~and~3.4.2, that the nodes $\{x_j\}_{1\le j \le n}$ of the Gaussian
quadrature of order $n$
are given by the zeros of the polynomial $p_n(x)$, and the weights are
given by
%
%e45 #&#
%
\begin{equation}
\label{weightsGaussianquadrature} w_j=\frac{a_n}{a_{n-1}}\frac{(p_{n-1},p_{n-1})_{\nu}}{p_{n-1}(x_j)p_n'(x_j)},
\end{equation}
where $a_k$ is the coefficient of $x^k$ in $p_k(x)$.

The following result demonstrates close connections between Gaussian
quadrature, orthogonal polynomials, Pad\'e approximations and Stieltjes
functions.

%th6 #&#
%
\begin{theorem}[(Theorems 2.2 and 3.1 in \cite{PadGau})]\label{PadeandGauss}
Consider a Stieltjes function
\[
f(z):=\int
_{[0,a]} \frac{\nu(\d x)}{1+xz}.
\]
Then
%
%e46 #&#
%
\begin{equation}
\label{n1npade} f^{[n-1/n]}(z)=\frac{(-z)^{n-1}q_{n-1}(-1/z)}{(-z)^np_n(-1/z)}=\sum
 _{i=1}^n
\frac{w_i}{1+x_i z},
\end{equation}
where $\{x_i\}_{1\le i \le n}$ and $\{w_i\}_{1\le i \le n}$ are the
nodes and weights of the Gaussian quadrature with respect to the
measure $\nu(\d x)$, $p_n(z)$ is the $n$th orthogonal polynomial with
respect to $\nu$ and $q_{n-1}(z)$ is the associated polynomial
of degree $n-1$, defined by
\[
q_{n-1}(z):=\int
_{[0,a]} \frac{p_n(z)-p_n(w)}{z-w} \nu(\d w).
\]
\end{theorem}

Next, we will discuss how one can compute the coefficients of the Pad\'
e approximation.
Consider a function $f(z)$ given\vspace*{1pt} by a formal series expansion
$f(z)=\sum_{i\ge0} c_i z^i$. Then the Pad\'e approximation
$f^{[m/n]}(z)=P_m(z)/Q_n(z)$ with $m\ge n$ can be found as follows
(provided it exists):
first, we solve the system of $n$ linear equations
%
%e47 #&#
%
\begin{equation}
\label{findb} \qquad \lleft[ %
\matrix{ c_{m-n+1} & c_{m-n+2}
& c_{m-n+3} & \cdots& c_{m}
\cr
c_{m-n+2} &
c_{m-n+3} & c_{m-n+4} &\cdots& c_{m+1}
\cr
c_{m-n+3} & c_{m-n+4} & c_{m-n+5} &\cdots&
c_{m+2}
\cr
\vdots& \vdots& \vdots& \ddots& \vdots
\cr
c_{m}
& c_{m+1} & c_{m+2} &\cdots& c_{m+n-1}} \rright]
\lleft[ %
\matrix{ b_n
\cr
b_{n-1}
\cr
b_{n-2}
\cr
\vdots
\cr
b_1} \rright] =- \lleft[ %
\matrix{ c_{m+1}
\cr
c_{m+2}
\cr
c_{m+3}
\cr
\vdots
\cr
c_{m+n}} \rright]
\end{equation}
and find $b_i$, $1\le i \le n$. These coefficients give us the
denominator $Q_n(z):=1+b_1z+b_2z^2+\cdots+b_nz^n$. Then, the
coefficients of the numerator
$P_m(z):=a_0+a_1z+a_2z^2+\cdots+a_mz^m$ can be calculated recursively
%
%e48 #&#
%
\begin{eqnarray}
\label{computecoefficientsa}
\nonumber
a_0&=&c_0,
\\
\nonumber
a_1&=&c_1+b_1c_0,
\\
a_2&=&c_2+b_1c_2+b_2c_0,
\\
\nonumber
& \vdots&
\\
\nonumber
a_m&=&c_m+\sum
_{i=1}^n
b_i c_{m-i}.
\end{eqnarray}
In practice, when $n$ is even moderately large, the system in (\ref
{findb}) will have a very large condition number, and solving the
system of linear equations (\ref{findb}) would involve a loss of
accuracy. This can be avoided by using higher precision arithmetic.
Another way to deal with this problem is to use expressions for Pad\'e
approximations given in terms of Gaussian quadrature [such as
(\ref{defpsin}), (\ref{phinkn}) and (\ref{phinkn2})]. There
exist several very fast and accurate methods for
computing the weights and nodes of the Gaussian quadrature; see \cite
{Gautschi,GolWel}.

Below we collect some other results on Pad\'e approximations, which are
used elsewhere in this paper.
%th7 #&#

\begin{theorem}[(Theorem 1.5.2 in \cite{Baker})]\label{thm152}
Given a formal series $f(z)=\sum_{i=0}^{\infty} c_i z^i$ and $a\neq
0$ we define $w=w(z)=az/(1+bz)$ and $g(w)=f(z)$. If the Pad\'e approximant
$f^{[n/n]}(z)$ exists, then $g^{[n/n]}(w)=f^{[n/n]}(z)$.
\end{theorem}

%th8 #&#
%
\begin{theorem}[(Theorem 1.5.3 in \cite{Baker})]\label{thm153}
Given a formal series $f(z)=\sum_{i=0}^{\infty} c_i z^i$ we define
$g(z)=(a+bf(z))/(c+df(z))$. If $c+df(0)\neq0$ and the Pad\'e approximant
$f^{[n/n]}$ exists, then
\[
g^{[n/n]}(z)=\frac{a+bf^{[n/n]}(z)}{c+df^{[n/n]}(z)}.
\]
\end{theorem}

%th9 #&#
%
\begin{theorem}[(Theorem 1.5.4 in \cite{Baker})]\label{thm154}
Assume that $k\ge1$ and $n$, $m$ are integers such that $n-k\ge m-1$.
Given a formal series $f(z)=\sum_{i=0}^{\infty} c_i z^i$ we define
\[
g(z)= \Biggl(f(z)-\sum
_{i=0}^{k-1}
c_i z^i \Biggr) z^{-k}.
\]
Then
\[
g^{[n-k/m]}(z)= \Biggl(f^{[n/m]}(z)-\sum
_{i=0}^{k-1}
c_i z^i \Biggr) z^{-k},
\]
provided either Pad\'e approximant exists.
\end{theorem}

%th10 #&#
%
\begin{theorem}[(Theorem 5.4.4 in \cite{Baker})]\label{thm544}
Let $f(z)$ be a Stieltjes series with radius of convergence $R>0$. Let
$A$ be a compact subset of
$\c\setminus(-\infty, -R]$. Define $\delta$ to be the distance
from $A$ to the set $(-\infty, -R]$ and $\rho:=R-\delta$. Then there
exists a constant $C=C(A)$ such that
for all $z \in A$ and all $n \ge1$ we have
\[
\bigl\llvert f(z)-f^{[n-1/n]}(z)\bigr\rrvert <C \biggl\llvert
\frac{\sqrt{\rho+z}-\sqrt{\rho
}}{\sqrt{\rho+z}+\sqrt{\rho}} \biggr\rrvert ^{2n}.
\]
\end{theorem}
\end{appendix}

\section*{Acknowledgments}
The authors would like to thank two anonymous referees for their
careful reading of the paper and for suggesting several improvements.

%\begin{supplement}[id=suppA]
%\slink[doi]{10.1214/00-AAPXXXXSUPP} %[doi,text={...}] - jei reikia
%suskaldyti doi
%\sfilename{aapXXXX\_supp.pdf}
%\end{supplement}

% imsref loaded by linak, 2015-02-06 10:02:24
%
% imsref loaded by linak, 2015-02-24 10:06:37

\printaddresses
\end{document}